\documentclass[twocolumn]{amsart}
\usepackage[utf8]{inputenc}
\usepackage{amsfonts, amssymb, amsmath, amsthm}
\usepackage{color}
\usepackage[normalem]{ulem}

\usepackage{epsfig,graphicx}
\usepackage[centering, includeheadfoot, hmargin=0.75in, tmargin=0.6in, 
  bmargin=0.75in, headheight=6pt]{geometry}

\newtheorem{theorem}{Theorem}
\newtheorem*{theorem*}{Theorem}

\newtheorem*{main-conjecture}{Main conjecture}

\theoremstyle{definition}

\def\DD{{\mathbb D}}  
\def\ZZ{{\mathbb Z}}

\def\RR{{\mathbb R}}

\def\NN{{\mathbb N}}

\def\eps{\epsilon}

\begin{document}

\title{From zero to positive entropy}
\author{Sylvain Crovisier}
\author{Enrique Pujals}
\maketitle

In the sciences in general, the phrase ``route to chaos'' has come to refer to a metaphor when some physical, biological, economic, or social system transitions from one exhibiting order to one displaying randomness (or chaos). Sometimes the goal is to understand which universal mechanisms explain that transition, and how one can describe systems that operate in a region between order and complete chaos. In other words, the goal is to understand the mathematical processes by which a system evolves from one
{\color{black} whose recurrent set is finite towards another one} exhibiting chaotic behavior as parameters governing the behavior of the system are varied. This has only been understood for one-dimensional dynamics. The present note exposes new approaches that allow one to move away from those limitations. 

A tentative global framework toward describing a large class of two-dimensional  dynamics,
inspired partially by the developments in the one-dimen\-sio\-nal theory of interval maps is discussed. More precisely, we present a class of intermediate smooth dynamics between one and higher dimensions. {\color{black} In this setting, it could be possible} to develop a similar one-dimensional type approach and {\color{black} in particular} to understand the transition from zero entropy to positive entropy.

\section*{Complexity in dynamics}
\begin{figure}[ht]
\begin{center}
\includegraphics[width=7cm,angle=0]{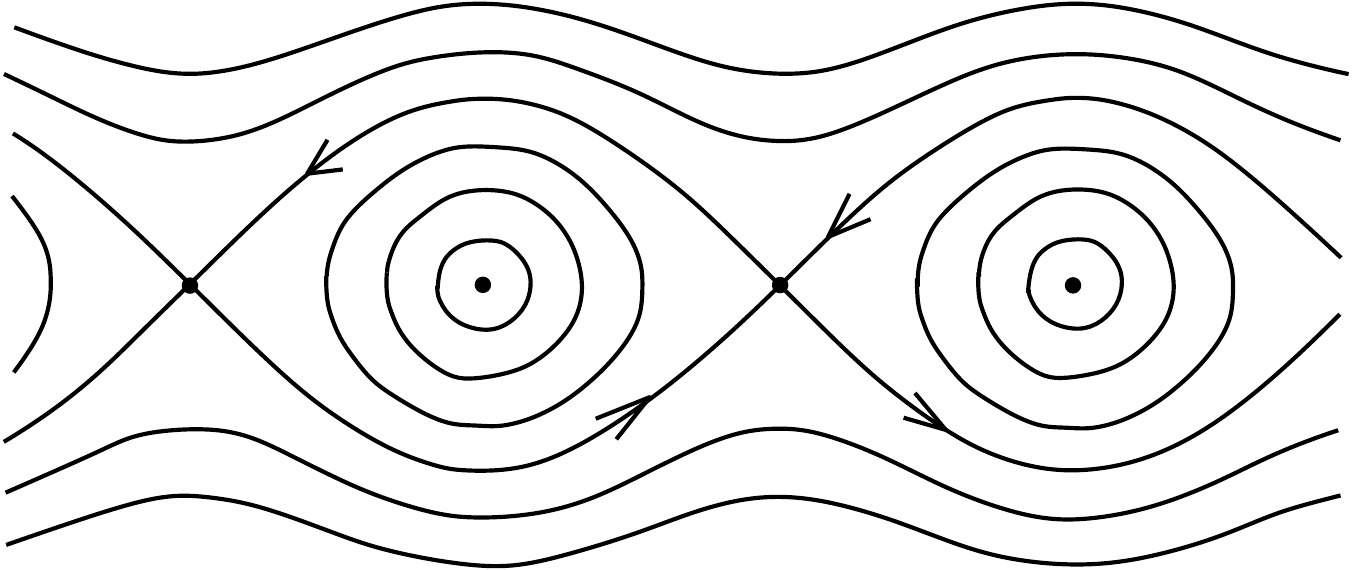}
\end{center}
\caption{A hamiltonian dynamics (the pendulum).\label{f.hamiltonian}}
\end{figure}

{\color{black} Considering a system which evolves in time,
the purpose of dynamical systems is} to describe the asymptotic behavior of its orbits.
As an example, one may {\color{black} think to} the gradient flow associated to a Morse function: there exists a finite number of equilibria and any other orbit is a curve which connects one equilibrium to another one.
\begin{figure}[ht]
\begin{center}
\includegraphics[width=4cm,angle=0]{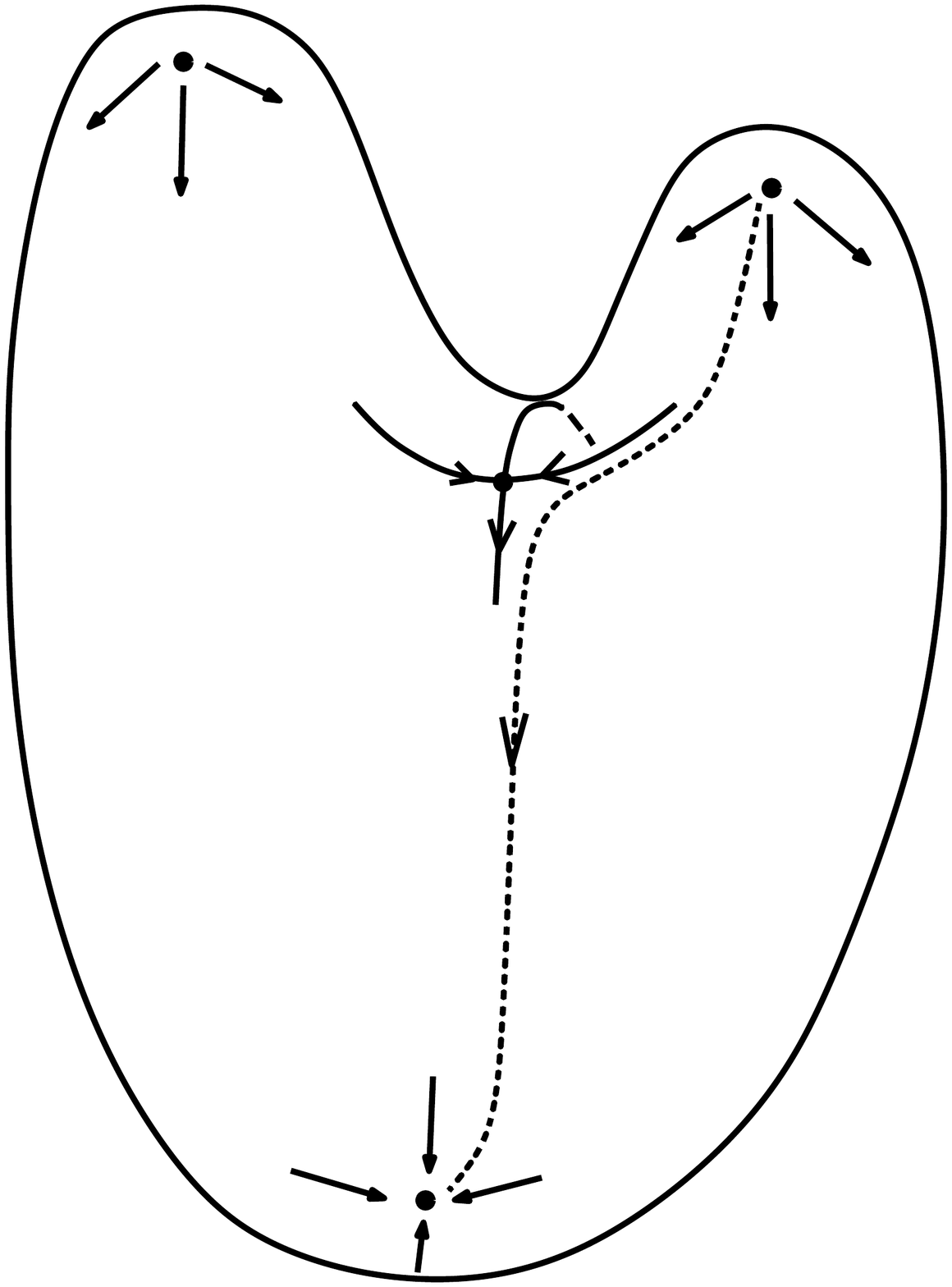}
\end{center}
\caption{A Morse-Smale dynamics defined by the gradient flow of a Morse function.\label{f.morse-smale}}
\end{figure}
One may also have in mind mechanical systems:
in the case of the ideal frictionless  pendulum, one gets a flow whose orbits are contained in the level sets of the energy function. See figures~\ref{f.morse-smale} and~\ref{f.hamiltonian}.

In this note we consider discrete time systems, defined by a map $f$ on a phase space $M$.
For instance, $f$ may be the time-$1$ map for the flows mentioned previously.
The forward orbits are the sequences of the form $x$, $f(x)$, $f(f(x))$, $f(f(f(x)))$,\dots
For convenience, one usually denotes $f^n(x)$ the image after $n$ compositions by $f$ and one of our goals is to characterize the accumulation sets of the orbits, usually called the limit sets.

\subsection*{The horseshoe map} 
\begin{figure}[ht]
\begin{center}
\includegraphics[width=5cm,angle=90]{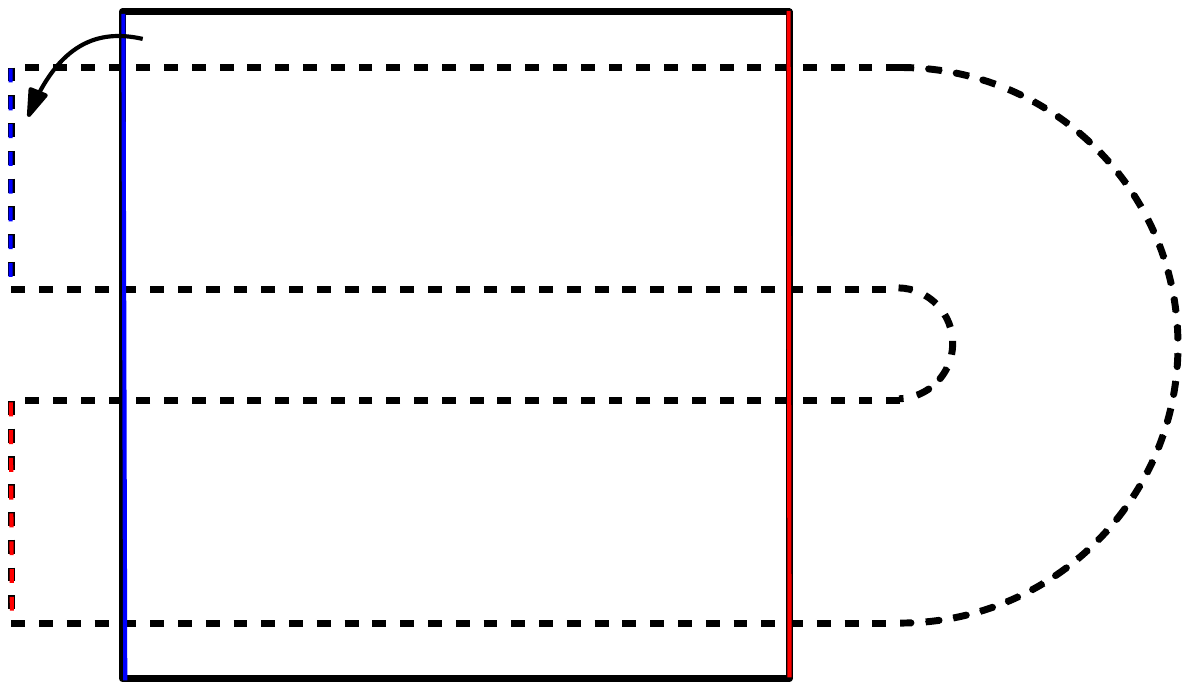}
\put(-84,1){$f$}
\put(-65,30){$R_0$}
\put(-30,30){$R_1$}
\put(-93,70){$R$}
\end{center}
\caption{The horseshoe map.\label{f.horseshoe}}
\end{figure}

For the above systems, {\color{black} or for others like rotations or isometries, the limit sets are very simple and the orbits are described easily.} But much richer behaviors exist.
This happens on surface, when a rectangle $R$ is vertically stret\-ched, horizontally contracted  and crossed
twice by its image $f(R)$. In this case $R\cap f(R)$ has two components ($R_0,R_1$) and any orbit contained in $R$ may be coded by a sequence in $\{0,1\}^\ZZ$, {\color{black} which represents the sequence of components met along the orbit.} See figure~\ref{f.horseshoe} and~\cite{what-is-horseshoe} for more details.
Conversely, any such sequence is realized by an orbit contained in $R$. This shows that for each time $n\geq 0$, at least $2^n$ different {\color{black} orbits} of the system may be distinguished at the scale of the rectangle $R$.

\subsection*{The topological entropy}
One measures the complexity of a dynamical system $f$ through its entropy. It is defined by fixing $\varepsilon>0$ and considering the
maximal number $N(\varepsilon,f,n)$ of {\color{black} orbits} under $f$ that can be distinguished at scale $\varepsilon$
{\color{black} up to} time $n$. The topological entropy $h(f)$ is the asymptotic exponential growth rate\footnote{\color{black} Note that $N(\varepsilon,f,n)$ increases as $\varepsilon$ gets smaller. Formally, one thus sets
$h(f)=\lim_{\varepsilon\to 0} \limsup_{n\to +\infty} \tfrac 1 n \log N(\varepsilon,f,n)$.}
of this quantity. It is always bounded when $f$ is a differentiable map of a compact manifold.
As we will see the dynamics differs dramatically when the entropy vanishes or is positive.

\subsection*{Morse-Smale dynamics}
For the simple systems pictured in figures~\ref{f.morse-smale} and~\ref{f.hamiltonian} the entropy is zero.
This is also the case for any Morse-Smale dynamics, i.e. for systems which generalize the gradient dynamics
in this way:
\begin{itemize}
\item there exist finitely many periodic orbits $O_1$,\dots, $O_\ell$,
each of them being hyperbolic (when $x\in O_i$ is fixed by $f^n$, $n\geq 1$,
the moduli of the eigenvalues of $Df^n(x)$ are different from $1$),
\item any other orbit accumulates in the past and in the future on two different orbits
$O_i,O_j$ with $i<j$.
\end{itemize}

\subsection*{Cascade of doubling periods and odometers}
An important example of a diffeomorphism on the disc with zero entropy {\color{black} and which is not Morse-Smale} has been first built in~\cite{GvST}
and exhibits infinite sequence of periodic orbits.
It can be described as follows: the disc $D$ is mapped into itself and is separated by a line
of points $\gamma$ whose forward orbit converge to a fixed point $x_0$.
The two components $D_l,D_r$ of $D\setminus \gamma$ are topological discs that are exchanged by the map
and contain points $x_l,x_r$ of a same $2$-periodic orbit.
Each disc $D_l$ or $D_r$ is divided by a line $\gamma_l$ or $\gamma_r$
of points whose orbit accumulates on $\{x_l,x_r\}$; the four components of
$D\setminus (\gamma_l\cup\gamma\cup \gamma_r)$ are cyclically permuted by the map
and each of them contains a point of a same $4$-periodic orbit.
The decomposition goes on inductively and produces one periodic orbit for each period $2^n$.

The collection of periodic points converges to an invariant Cantor set $K$. The restriction of the dynamics
to $K$ is conjugated to the addition by $1$ on the group of dyadic integer $\mathbb{Z}_2$,
and for that reason the limit Cantor set is called an odometer (or adding machine).
\begin{figure}[ht]
\begin{center}
\includegraphics[width=6cm,angle=0]{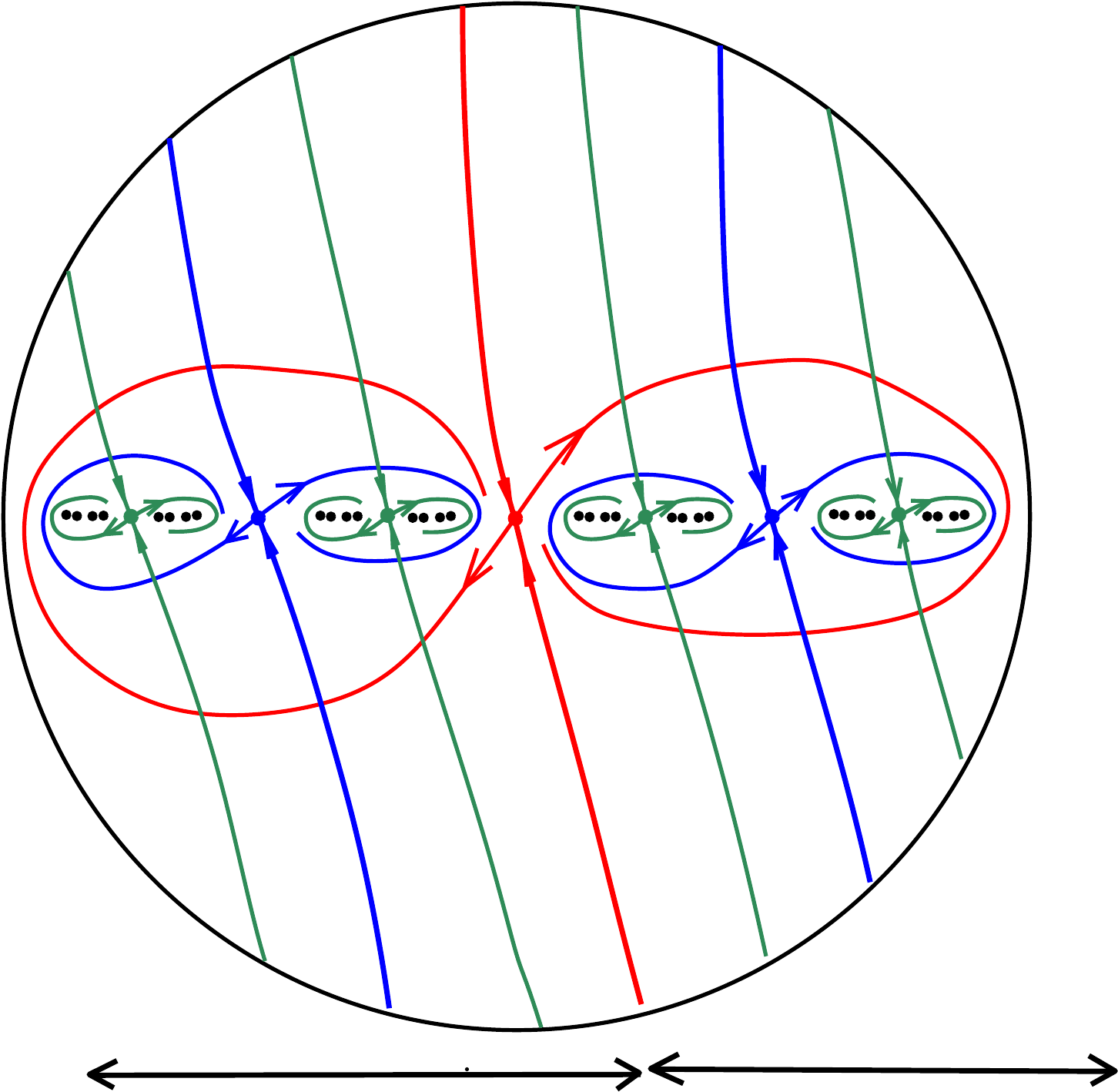}
\put(-93,95){\color{red}$x_0$}
\put(-88,40){\color{red} $\gamma$}
\put(-43,55){\color{black} $\gamma_r$}
\put(-118,50){\color{black} $\gamma_l$}
\put(-53,96){\color{black}$x_r$}
\put(-130,97){\color{black}$x_l$}
\put(-130,6){\color{black}$D_l$}
\put(-40,6){\color{black}$D_r$}
\put(-10,70){$D$}
\end{center}
\caption{A cascade of doubling periods accumulating on an odometer.\label{f.odometer}}
\end{figure}
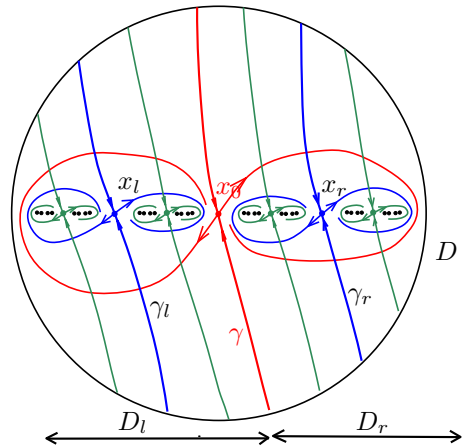

\subsection*{Transverse homoclinic intersections}
\begin{figure}[ht]
\begin{center}
\includegraphics[width=6cm,angle=0]{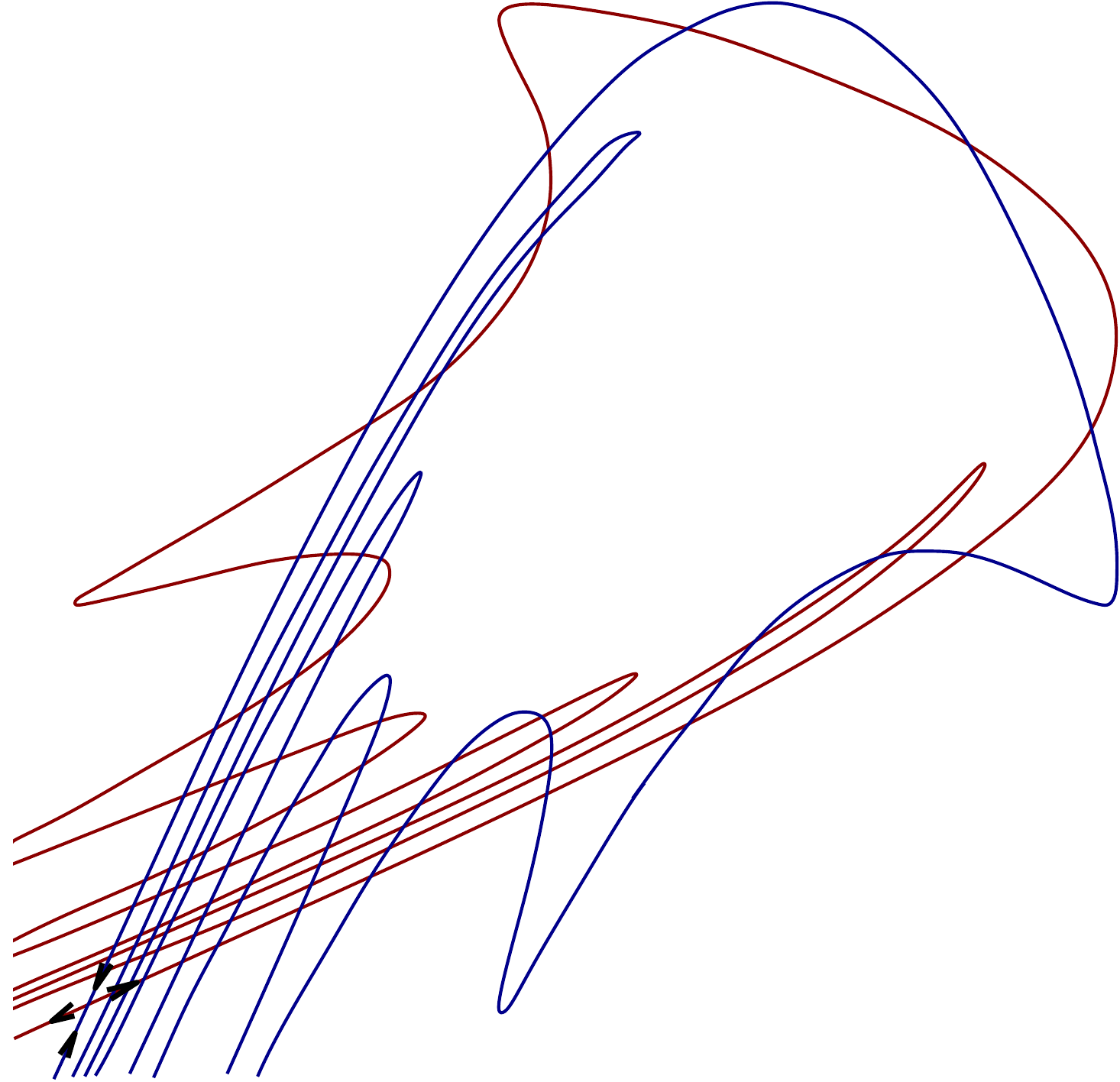}
\end{center}
\caption{The net formed by the stable and unstable manifolds of a fixed point, when there exists a transverse homoclinic intersection.
\label{f.homoclinic}}
\end{figure}

For the horseshoe map, pictured on figure~\ref{f.horseshoe}, we have seen that the number of iti\-nera\-ries {\color{black} at time $n$ has the lower bound $2^n$} and the entropy is at least $\log(2)$. A generalization of that phenomenon occurs frequently in differentiable dynamics.
Indeed, let us consider a diffeomorphism with a {\color{black} fixed point $p$
which is a \emph{hyperbolic saddle}: the tangent space at $p$ decomposes as the sum of two invariant
subbundles $T_pM=E^s\oplus E^u$, such that the eigenvalues of $Df|_{E^s}$ (resp. $Df|_{E^u}$)
have a modulus smaller than $1$ (resp. larger than $1$). Then the set of points whose forward (resp. backward) orbit converges to $p$
is an immersed  submanifold $W^s(p)$ (resp. $W^u(p)$), called stable manifold (resp. unstable manifold) of $p$.
Poincar\'e has noticed a fascinating phenomenon: when these manifolds have a transverse intersection point (different from $p$ itself),
then they have to intersect in an intricate way, see figure~\ref{f.homoclinic}.}

Smale has proved that such a transverse homoclinic intersection forces $f$ to have a horseshoe, implying that the entropy is positive.
{\color{black} For periodic orbits $O$ with period $\tau$ which are saddle (i.e. which split into saddle fixed points of $f^\tau$), one defines analogously the stable and unstable manifolds $W^s(O), W^u(O)$;
a transverse intersection (outside $O$) between them implies that $f$ admits a horseshoe.}
Katok has shown the converse in dimension~$2$. {\color{black} On surfaces, positive entropy is thus equivalent
to the existence of a horseshoe:}
\medskip

\noindent
{\bf Theorem} (Smale, Katok).
\emph{For $C^2$-diffeo\-mor\-phisms on a surface, the topological entropy is positive if and only if
there exists a hyperbolic periodic orbit with a transverse homoclinic point.}
\smallskip

{\color{black} In higher dimensions, only Smale's implication remains,
but one may ask if Katok's one still holds for ``most" diffeomorphisms.}

\section*{Transition to chaos in the space of systems}
It appears that the space of diffeomorphisms splits into two classes with very different dynamics:
those with zero entropy and those with positive entropy.
Each of these classes contains open sets: the set of Morse-Smale diffeomorphisms on the one hand,
and the set of systems exhibiting a transverse homoclinic orbit on the other hand.
This naturally raises the following questions.
\medskip

\noindent
\textbf{Q1.}
\emph{Is the set of Morse-Smale diffeomorphisms dense in the set of systems with zero entropy?}
\medskip

\noindent
\textbf{Q2.}
\emph{Is the set of diffeomorphisms with a transverse homoclinic intersection dense in the set of systems with positive entropy?}
\medskip

If these questions have positive answers, the interface of these two classes is small and
one goal would be to understand the systems at the transition between simple {\color{black} (zero entropy)} and complicated {\color{black} (positive entropy)} dynamics. In particular we would like to identify, if it exists, the phenomena that generates entropy.
\medskip

\noindent
\textbf{Q3.}
\emph{Can one characterize systems that belong to the boundary of the class of dynamics with zero entropy?}
\medskip

We note that the two open classes introduced before can be distinguished by the number of periodic orbits present in the system:
it is stably finite in one case, and stable infinite in the other case.
\medskip

\noindent
\textbf{Q4.}
\emph{Can one identify the transition from finitely to infinitely periodic orbits?}
\medskip

We will discuss these questions in different settings, starting with the lower dimensions.

\section*{One-dimensional  dynamics}
These questions have  already been addressed in dimension $1$.
We will focus on continuous maps acting on the {\color{black} closed} interval, but one could also consider maps acting on other one-dimensional spaces like the circle or trees. For monotone maps, it is not difficult to prove that the accumulation points of any orbit is either given by fixed points
(if the map is increasing) or by a unique fixed point and periodic orbits of period two (if the map is decreasing).

A richer situation holds with non-invertible maps. The action of quadratic polynomials on the real line illustrates how different possible dynamical scenarios arise.
Without loss of generality, one can consider the quadratic family $$\color{black} f_a\colon x\mapsto ax(1-x)$$ which satisfies $f_a([0,1])\subset [0,1]$ for $a\in [0,4]$. There is a value $a_\ast$ such that the topological entropy vanishes when $a\leq a_\ast$ and is positive for the other parameters.
{\color{black} As $a<a_\ast$ increases the (finite) number of periodic orbits increases and $a_\ast$ is the smallest parameter exhibiting periodic points
with arbitrarily large periods.}

The natural ordering of the interval allows a combinatorial approach.
For instance, exploring the richness of that total order structure,
Milnor and Thurston have developed a Kneading Theory, giving a complete
description of all topological possibilities for the dynamics of a family of endomorphisms of the interval,
with a given number of monotonicity branches.

\subsection*{Periodic approximation in the interval}

{\color{black} A point $x_0$ is recurrent if its forward orbit meets any of its neighborhoods; this is the case when $x_0$ is periodic.}
A result, that highlights the strength of the order structure, asserts:
\medskip

\noindent
\textbf{Property} (L.-S. Young~\cite{LSY}).
\emph{For interval maps, the periodic points are dense in the recurrent set.}
\medskip

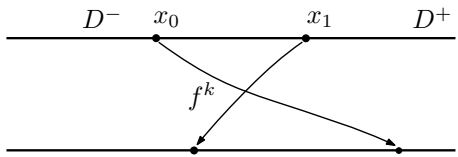
\begin{figure}[ht]
\begin{center}
\includegraphics[width=6cm,angle=0]{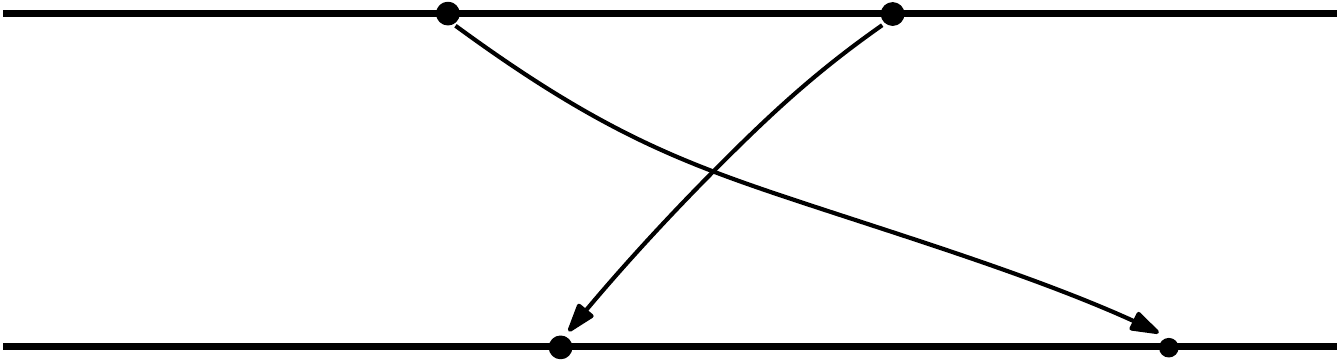}
\put(-105,19){ $f^k$}
\put(-60,50){ $x_1$}
\put(-118,50){ $x_0$}
\put(-145,48){ $D^-$}
\put(-20,48){ $D^+$}
\end{center}
\caption{Localization of a periodic point in the interval.
\label{f.closing-interval}}
\end{figure}

This fact is unknown in higher dimension, even from a smooth generic point of view.  A simple proof of this fact goes along the following lines.  Let us consider {\color{black} a recurrent (non-periodic) point $x_0$ and a forward iterate $x_1 = f^n(x_0)$ close to $x_0$.} Without loss of generality, one can assume that $x_0 < x_1$. We claim that there is a periodic point for $g := f^n$
inside $(x_0, x_1)$.
{\color{black} Indeed one can easily check that $x_0$ is still recurrent for $g$.
Since $x_0 < x_1$, there exists a positive integer $k$ such that $g^k(x_1) < x_1$.
Taking the smallest $k$ also gives} $g^k(x_0) = g^{k-1}(x_1) > x_1.$ Therefore we have a continuous map $g^k: [x_0, x_1]\to [0, 1]$ such that $g^k(x_0) > x_1$ and $g^k(x_1) < x_1.$ Hence the graph of $g^k$ crosses the diagonal inside $(x_0, x_1)$: there is a point $p\in (x_0, x_1)$ which is fixed for $g^k$ and so periodic for $f$, see figure~\ref{f.closing-interval}.

We can recast the previous proof, avoiding an explicit use of the order, and in a way that can be generalized to other contexts. Under the same choices of $x_0$ and $x_1$ as above, we define the intervals $D^-=[0, x_0]$ and $D^+=[x_1, 1]$ and we take the first positive integer $k$ such that $g^k(x_1)\notin D^+.$
Such an integer exists since $x_0$ is recurrent and does not belong to $D^+$.  By the choice of $k$ observe that $g^k(x_0)= g^{k-1}(x_1) \geq x_0$
Let $h\colon [0, 1]\to [x_0, x_1]$ be a continuous map which coincides with the identity on $[x_0,x_1]$
and such that $h([0,x_0])=x_0$, $h([x_1,1])=x_1$.
Then the map $h\circ g^k:[x_0, x_1]\to [x_0,x_1]$ has a fixed point $p\in [x_0, x_1].$
Note that $h\circ g^k(x_0)=x_1$ (since $g^k(x_0)\in D^+$) and $h\circ g^k(x_1)\neq x_1$ (since $g^k(x_1)\not\in D^+$).
Therefore $p$ belongs to $(x_0, x_1)$. Since $h$ is the identity on $(x_0, x_1)$, the point $p$ is a fixed point of $g^k$.

\subsection*{One-dimensional dynamics and zero entropy}
A simple characterization of positive entropy in the interval has been given by Misiurewicz:
\medskip

\noindent
\textbf{Property} (Misiurewicz).
\emph{An interval map has positive entropy if and only if there exists two disjoint intervals such that the  images  of each interval
by an iterate $f^\ell$ of the map  contains the union of both interval.}
\medskip

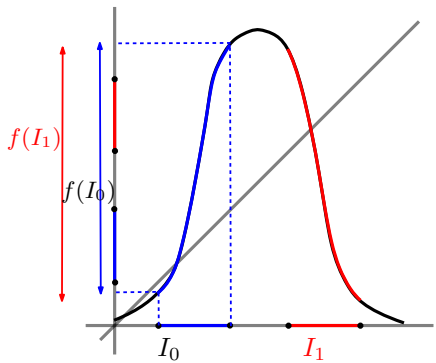
\begin{figure}[ht]
\begin{center}
\includegraphics[width=5cm,angle=0]{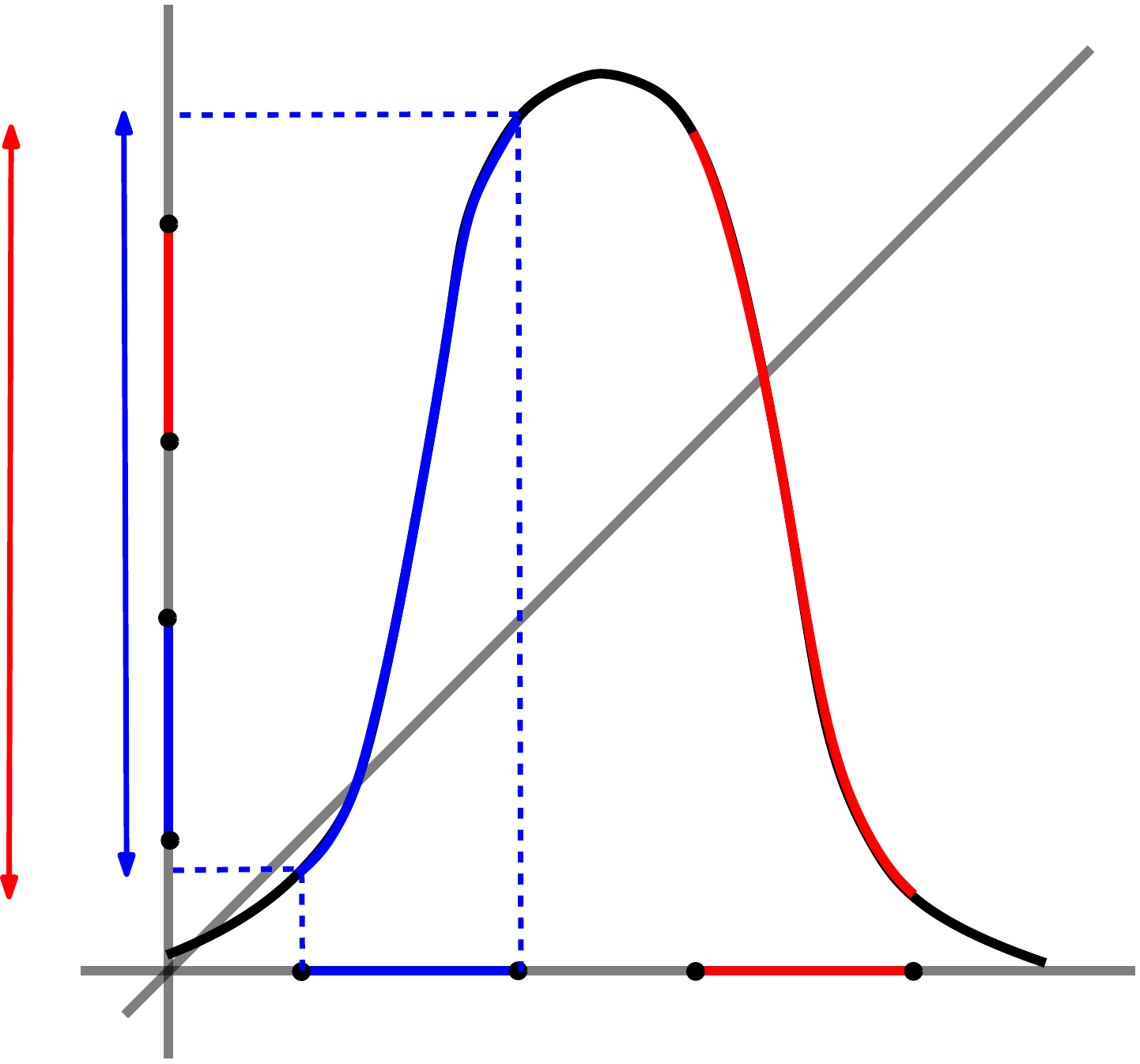}
\put(-105,0){\color{black} $I_0$}
\put(-50,0){\color{red} $I_1$}
\put(-162,80){\color{red} \small $f(I_1)$}
\put(-141.5,60){\color{black} \small $f(I_0)$}
\end{center}
\caption{Map satisfying Misiurewicz criterion.
\label{f.misiurewicz}}
\end{figure}

The reason is analogous to {\color{black} Smale-Katok's theorem}:
the number of itineraries for $f^\ell$ with respect to these intervals grows as $2^n$, see figure~\ref{f.misiurewicz}.

Another historical result in this combinatorial theory is Sarkovskii's hierarchy of periodic orbits. It implies:
\medskip

\noindent
\textbf{Property} (Sharkovskii).
\emph{Interval maps with zero entropy only admit periodic points of period $2^n$, $n\geq 0$.}
\medskip

Let us discuss {\color{black} Sharkovskii's property} in the case of \emph{unimodal maps} of the unit interval $[0,1]$,
i.e. continuous maps with only one turning point $c\in [0,1]$, one strictly increasing interval $[0,c]$
and one strictly decreasing interval $[c,1]$. The following dichotomy then holds:
\medskip

\noindent
\textbf{Property}.
\emph{For unimodal maps $f$ with zero entropy,\\
-- either all forward {\color{black} orbits} converge to a fixed point,\\
-- or  $f$ is {\em renormalizable}: there exists an interval $I$ containing $c$ such that $f(I)\cap I=\emptyset$, $f^2(I)\subset I$, $f^2_{|I}$ is unimodal and the forward orbit of any point either converges to a fixed point or enters in the interval $I$.}
\medskip

In particular, any periodic orbit is either fixed, or has even period. Applying the property to $f^2_{|I}$ shows that $4$ divides any period larger than $2$. Arguing inductively, one concludes that the allowed periods have the form $2^n$.
 
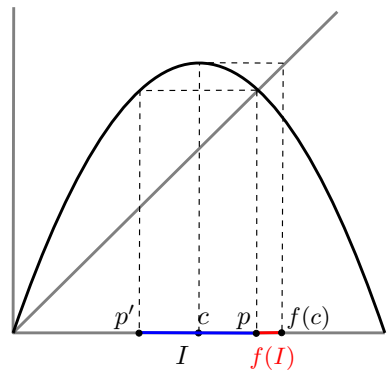
\begin{figure}[ht]
\begin{center}
\includegraphics[width=5cm,angle=0]{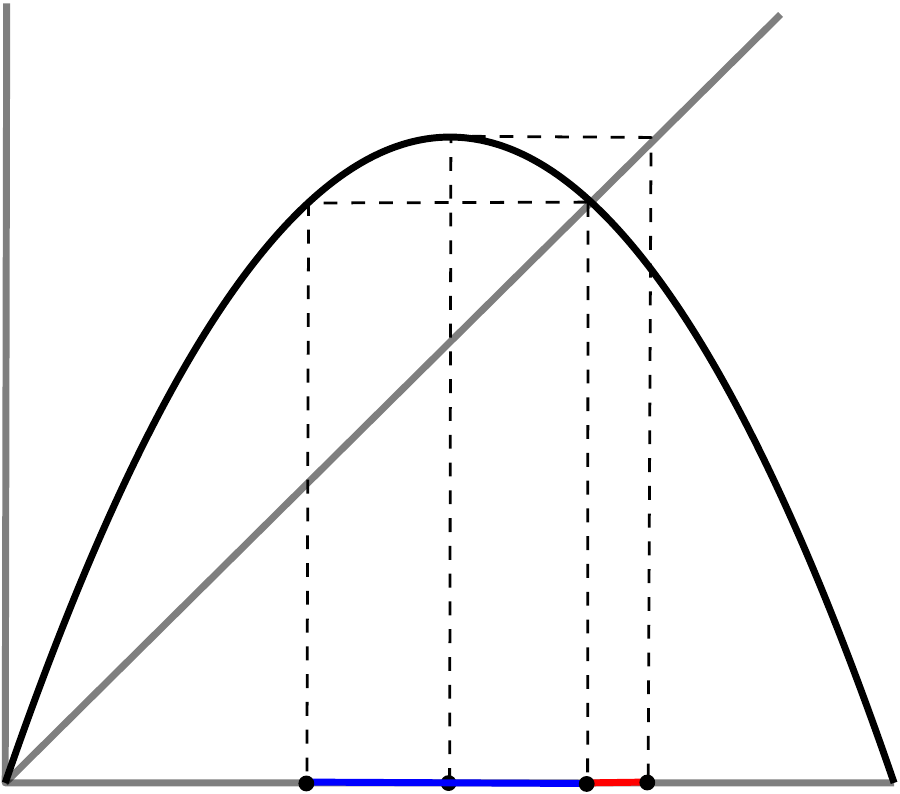}
\put(-80,-10){\color{black} $I$}
\put(-52,-10){\color{red} \small $f(I)$}
\put(-75,5){ $c$}
\put(-60,5){ $p$}
\put(-38,5){\small $f(c)$}
\put(-103,5){$p'$}
\end{center}
\caption{Unimodal map which is renormalizable.
\label{f.unimodal}}
\end{figure}

 The dichotomy can be obtained by considering separately the two cases $f(c)\leq c$ and $f(c)>c$. In the first case, $f([0,c])\subset [0,c]$ and since $f|_{[0,c]}$ is increasing, it follows that any orbit in $[0,c]$ converges to a fixed point; since $f([0,1])\subset [0,c]$,
 the same property holds on the whole interval $[0,1]$. In the second case, observe that
 since $f(c)>c$ and $f(1)\leq 1$, there is a fixed point $p\in (c,1]$. Let us introduce the maximal interval $I:=(p',p)\subset (0,p)$
 whose image is contained in $(p,1)$ (note that either $p'=0$ or $f(p')=p$). Observe that the turning point $c$ belongs to that interval, $f(I)\cap I=\emptyset$ and $f^2|_I$ is unimodal. Also $f^2(p')\leq f^2(p)=p$. Since the entropy is zero, it follows from Misiurewicz property that $f^2(c)\in I$ and therefore $f^2(I)\subset I.$ It remains to see that any forward orbit either converges to a fixed point or enters inside $I$. Note first that if $x<p_{-}$ then $f(x)< f(p')=p $ and so either $f(x)\geq p'$ (in this case $f(x)$ is contained in the interval $I$) or $f(x)< p'$ and the argument can be repeated: if the forward orbit of $x$ does not enter in $I$, it remains in the increasing interval $[0,c]$ and converges to a fixed point. In the last case $x>p$: the image $f(x)$ belongs to the increasing part and we are reduced to the first case. See  figure~\ref{f.unimodal}.
 
\subsection*{Infinite renormalization and odometers.}
From the previous discussion, one concludes that for unimodal maps with zero entropy,
two cases are possible.

A first possibility is that the inductive renormalization described in the previous paragraph
stops after a finite number $m$ of steps. The set of periods is then the
finite set $\{2^k,\; 0\leq k\leq m\}$. Any forward orbit accumulates on one periodic orbit.
The dynamics is similar to a Morse-Smale dynamics (although the number of periodic points
of a given period may be infinite).

Otherwise one says that $f$ is \emph{infinitely renormalizable}.
For each $k\geq 0$, let $I_k$ denote a renormalization interval with period $2^k$, so that
$V_k:=I_k\cup f(I_k)\cup\dots\cup f^{2^k-1}(I_k)$ is forward invariant.
The family $(V_k)$ is decreasing and the intersection
is an invariant compact set $\mathcal{K}$. When $f$ is $C^2$ and $D^2f(c)\neq 0$,
it is a Cantor set\footnote{It follows from the \emph{no wandering interval theorem}, see~\cite{dMvS}.},
and the dynamics on $\mathcal{K}$ is the same as in the example of  figure~\ref{f.odometer}: it is an odometer.
Inside such a set, all the orbits are dense and follows the same statistic: they distribute toward the same invariant probability measure, and visit a set $I_k$ with frequency $2^{-k}$. Any forward orbit of $f$ accumulates either on a periodic orbit, or on the odometer.

\subsection*{The renormalization operator.}
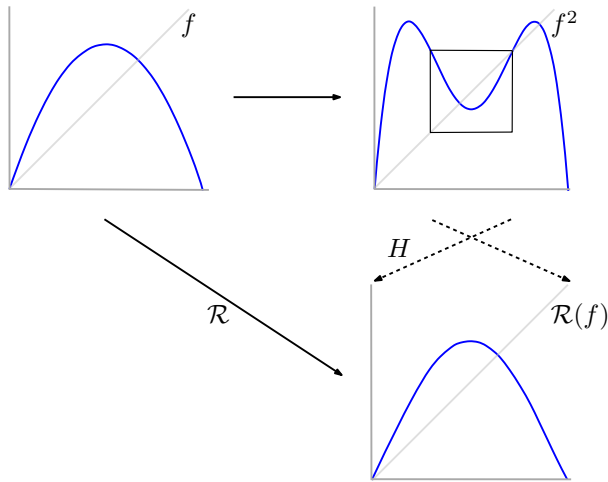
\begin{figure}[ht]
\begin{center}
\includegraphics[width=7.5cm]{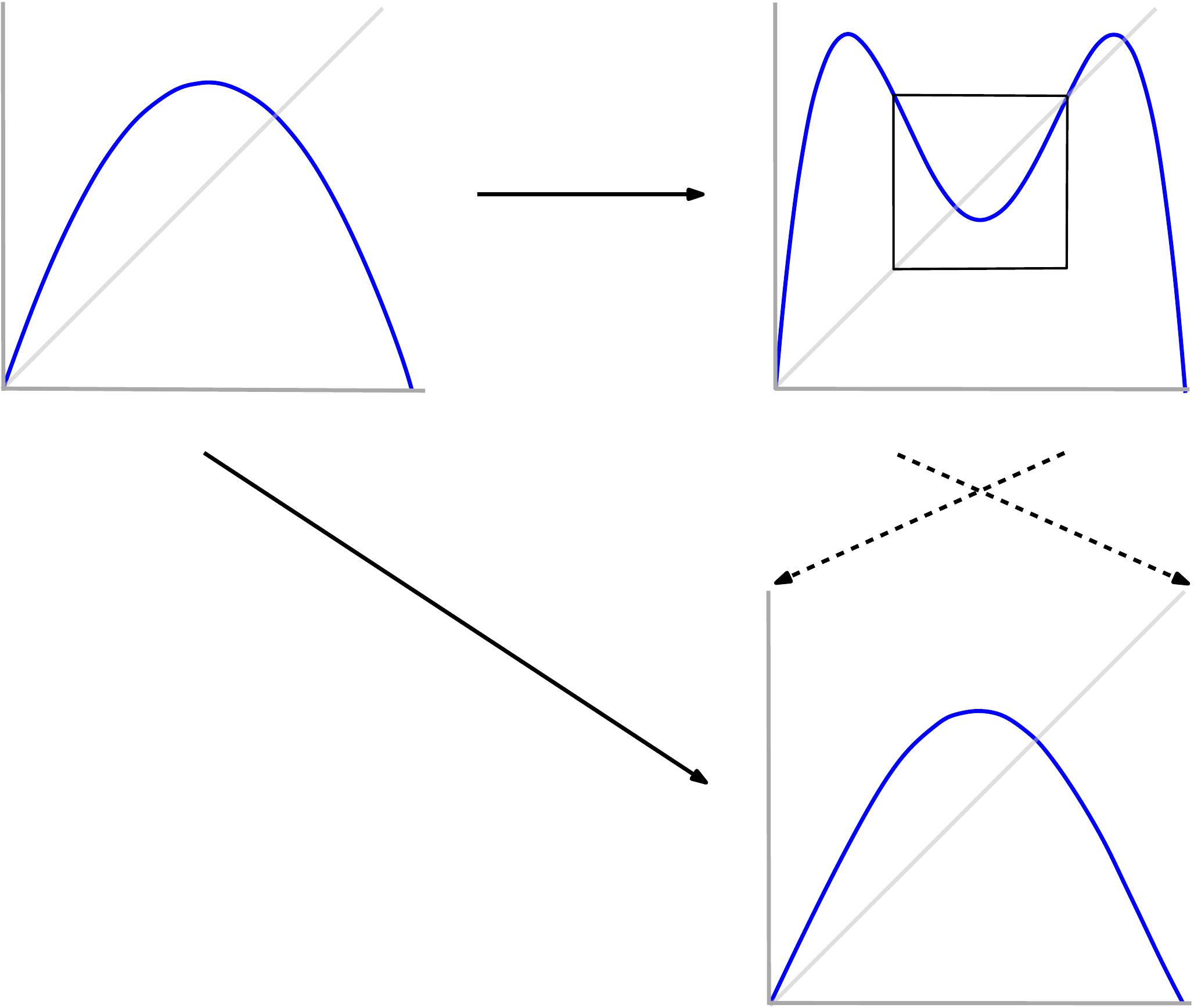}
\put(-138,60){$\mathcal{R}$}
\put(-70,85){$H$}
\put(-148,170){$f$}
\put(-8,170){$f^2$}
\put(-8,60){$\mathcal{R}(f)$}
\end{center}
\caption{Renormalization of a unimodal map.
\label{f.renormalization}}
\end{figure}

Deepening the idea of renormalization, Coullet-Tresser and independently Feigenbaum, have proposed to consider the renormalization operator $\mathcal{R}$ acting on the space of {\color{black} smooth} unimodal maps {\color{black} with a quadratic turning point}: to any map which is renormalizable on a maximal interval $I$,
it associates the map $\mathcal{R}(f):=H\circ f^2|_I\circ H^{-1}$, where $H$ is the orientation-reversing affine homeomorphism between $I$ and $[0,1]$. See figure~\ref{f.renormalization}.

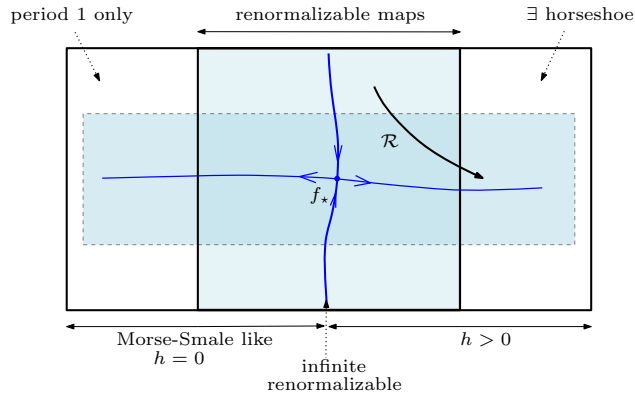
\begin{figure}[ht]
\begin{center}
\includegraphics[width=7cm,angle=0]{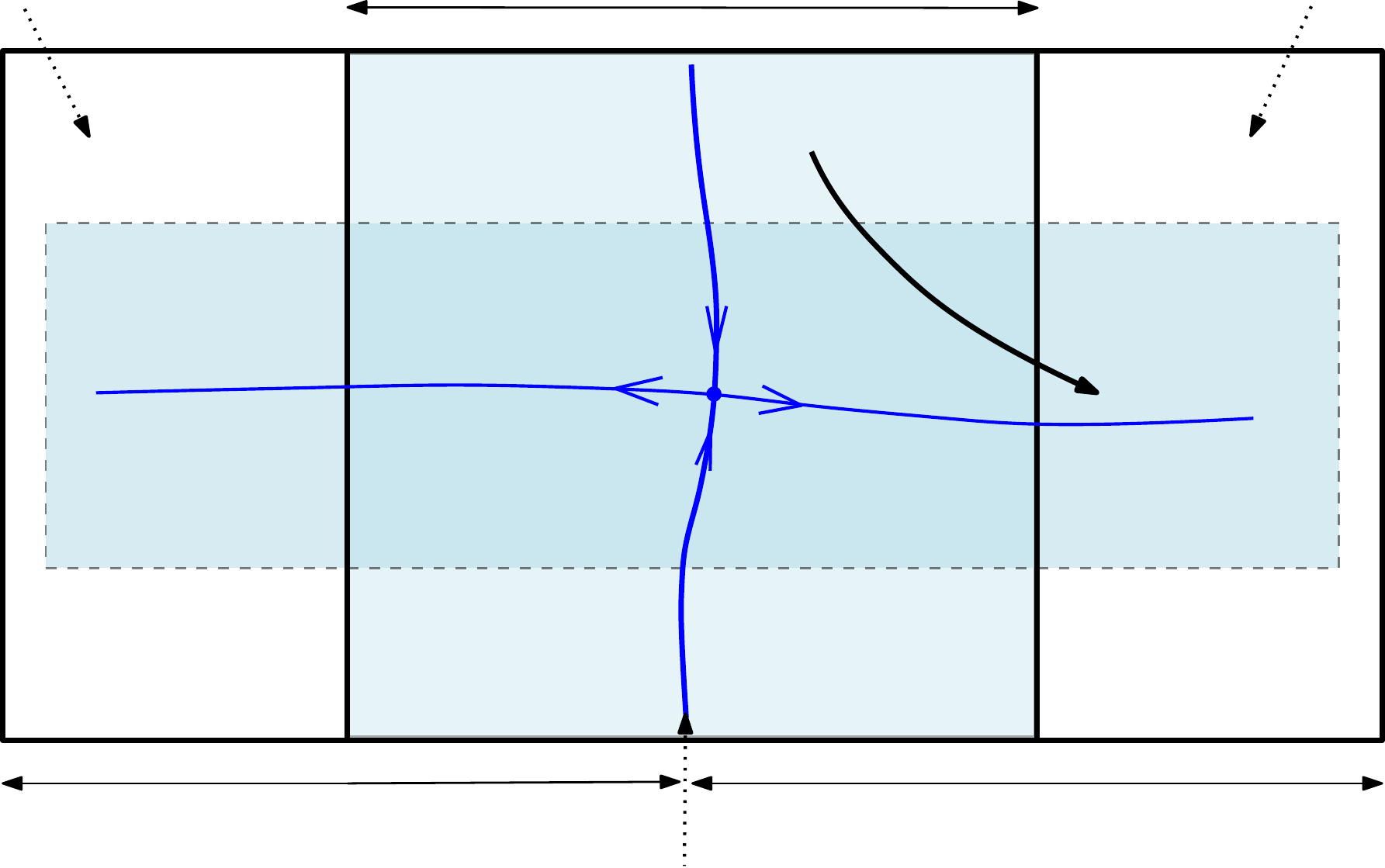}
\put(-110,-5){\scriptsize infinite}
\put(-123,-12){\scriptsize renormalizable}
\put(-180,5){\scriptsize Morse-Smale like}
\put(-166,-2){\scriptsize $h=0$}
\put(-50,5){\scriptsize $h>0$}
\put(-107,60){\color{black} \scriptsize $f_\star$}
\put(-80,80){\scriptsize $\mathcal{R}$}
\put(-135,128){\scriptsize renormalizable maps}
\put(-25,128){\scriptsize $\exists$ horseshoe}
\put(-220,128){\scriptsize period $1$ only}
\end{center}
\caption{Dynamics of the renormalization operator $\mathcal{R}$ on the space of unimodal maps with quadratic turning point.
\label{f.renormalization-global}}
\end{figure}

{\color{black} These people have realized} that the dynamics of $\mathcal{R}$ is the key for understanding the boundary of the set of maps with zero entropy.
They conjectured that the operator
has a unique fixed point $f_\star$, which is hyperbolic.
The set of unimodal maps whose sequence of renormalization {\color{black} converges} to $f_\star$ is a one-codimensional submanifold (which correspond to infinitely renormalizable maps). Outside this stable manifold,
the renormalizations stop after finitely many steps. On one side the dynamics has a Morse-Smale behavior: the number of periods is finite and the entropy vanishes. On the other side of the stable manifold, the dynamics renormalizes until a horseshoe appears and
the entropy is positive. See figure~\ref{f.renormalization-global}.

The definite mathematical proof of these results started {\color{black} first in the analytic context} with Sullivan's program \cite{Su}, approaching the Feigenbaum-Coullet-Tresser Renormalization Conjecture based on Teichm\"uller theory, and finished with the proof by Lyubich \cite{L}, showing the hyperbolicity of the renormalization fixed point; this has been later extended to lower regularity in~\cite{dFdMP}.  Partial results about maps with more monotonicity branches (multimodal maps) and the associated transition to chaos have been obtained  by many authors (see \emph{e.g.} \cite{MiTr} and references cited or citing).

These results also explains some quantitative and universal phenomena appearing {\color{black} when the system changes inside one-parameter families.} Every family of unimodal maps presents essentially the same dynamical features as it passes from {\color{black} zero} to {\color{black} positive entropy: for instance when one measures the size of the set of parameters for which some periods appear.}
This is sometimes called topological universality for one-dimensional dynamics {\color{black} since} it allows to show that the quadratic family encapsulates all possible dynamical behaviors.

\section*{Dissipative surface dynamics}
One can naively wonder if that rich and meaningful description of the dynamics on the interval, can be extended to higher dimension. The next level of complexity to be considered are dissipative two-dimensional invertible maps acting on the disc, i.e. diffeomorphisms $f$ from the $2$-disc $\DD$ into its image $f(\DD)\subset \DD$ and which contract the volume. Therefore, the iterates of disc are confined to a set whose $2$-dimensional volume vanishes and which seems to have a one-dimensional structure.

However,  there are phenomena in this setting with no one-dimensional counterpart: there exists a residual\footnote{The residual sets refer here to the Baire category: the phenomena holds on a G$_\delta$ set which is dense inside a non-empty open set of $C^2$-diffeomorphisms.} set of dissipative diffeomorphisms of the disc exhibiting infinitely many attracting periodic orbits with arbitrarily large periods (this property is called the Newhouse phenomenon, see \cite{N}), whereas generic smooth one-dimensional maps have only finitely many attracting periodic points\footnote{It follows from the generic finiteness of attractors, see theorem~B' in~\cite{dMvS}}.

\subsection*{H\'enon maps}
\begin{figure}[ht]
\begin{center}
\includegraphics[width=5cm]{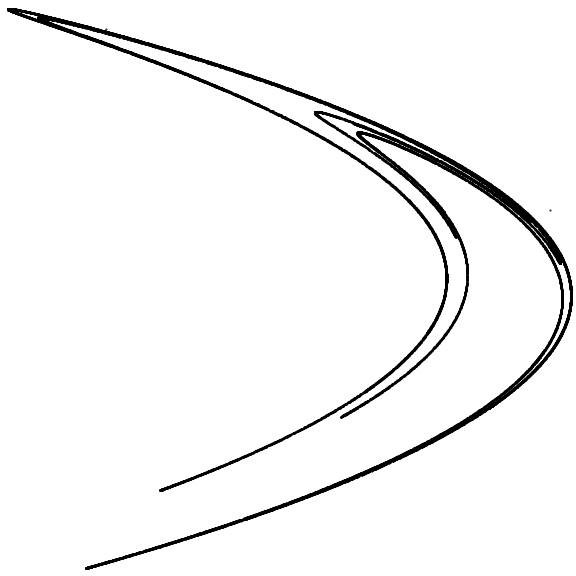}
\end{center}
\caption{The forward orbit of the point $(0.35,0.35)$ under the H\'enon map $(x,y)\mapsto(1-1.4 x^2+y, 0.3 x)$.
\label{f.henon}}
\end{figure}

One classical example of  dissipative surface maps is the H\'enon map which is defined by the formula $(x,y)\mapsto(1-ax^2+y, by)$ where $a$ and $b$ are real parameters and $b$ has modulus in $(0,1)$. See figure~\ref{f.henon}.

This family was introduced by H\'enon back in the seventies as a non-linear model {\color{black} displaying} complicated dynamics. In the age of computers and of computer graphics, H\'enon maps are one of the simplest two-di\-men\-sio\-nal systems used, through numerical simulations, to show how iterations produce extraordinarily complex behaviors.  However, the phenomena observed computationally have been rigorously explained
{\color{black} only for very small sets of parameters.}

Observe that $b$ is the Jacobian of the map and {\color{black} by letting $b$ tend to zero, one} recovers the classical quadratic family described.
One may thus expect that the two-dimensional maps have dynamical features of the interval quadratic map $x\mapsto 1-ax^2$,
even if the H\'enon maps also display new properties, such as the intriguing Newhouse's phenomenon.

\subsection*{Difficulties}
In practice, the two-dimensional systems are much more difficult to describe and much less is known for the H\'enon maps than for the quadratic family.
One reason is that there does not exist an obvious ordering on the phase space as in dimension $1$, so that a combinatorial structure of the dynamics is much more difficult to introduce. In particular there is no point in the disc which generalizes a priori the turning points in the interval.
A notion of critical points may be defined for some parameters but their number turns out to often be infinite while there is only one for quadratic maps.  

In that sense, trying to look for a general description of dissipative diffeomorphisms on the disc (with the necessary adaptations) as it has been performed in the one-dimensional context, could be considered as overambitious and unattainable with such a level of generality, or  even looking to the wrong paradigm.

\subsection*{Perturbative approaches {\color{black} and} strong dissipation.}
Through deep analysis it is possible to describe subsets of parameters inside the H\'enon family as small perturbations of the one-dimensional setting, either with positive entropy~\cite{BC} or with zero entropy~\cite{dCLM}. These works necessarily suppose {\color{black} that} the Jacobian $b$ is extremely close to $0$. {\color{black} This setting} will be qualified as ``strongly dissipative regime".

\subsection*{Zero entropy -- a conjecture by Tresser}
The Newhouse examples mentioned above have {\color{black} transverse homoclinic orbits and} positive entropy, so it is possible that 
when the entropy vanishes, the differences between interval dynamics and dissipative dynamics in the disc may disappear.
This expectation is encapsulated in a conjecture by Tresser~\cite{GT}.
It focuses on maps at the bifurcation locus between zero and positive entropy.
The natural examples are the diffeomorphisms pictured in figure~\ref{f.odometer}:
similarly to the one-dimensional case,
for each positive integer $n$ there is exactly one periodic orbit with period $2^n$, and no other period {\color{black} exists}.
\medskip

\noindent
{\bf Conjecture} (Tresser).
{\em In the space of dissipative diffeomorphisms of the disc, generically,
maps which belong to the boundary of the {\color{black} subset} of systems with zero entropy have an infinite set of periodic orbits
with periods $m.2^k$, for a given $m\geq 1$ and all $k \geq 0.$}
\medskip

In other terms, it asserts that at the transition between zero and positive entropy,
there exists a doubling cascade of periodic orbits.

\section*{Mildly dissipative surface dynamics}
After recognizing these {\color{black} difficulties of the dissipative surface dynamics}, we now present an open large class of dissipative diffeomorphisms acting on the disc $\mathbb{D}$, that has been introduced in~\cite{CP} and that
that captures {\color{black} key} properties of one-dimensional maps: abundance of periodic points in the recurrent set, order structure through one-dimensional reduction, renormalization structure in the entropy zero case, etc. However, it keeps two-dimensional features, showing all the well known complexity of dissipative surface diffeomorphisms; moreover it includes the H\'enon family with Jacobian $b$ up to $1/4$ and therefore goes beyond classic perturbative strategies.

\subsection*{Mild dissipation}
As mentioned before, the theory of real one-dimensional dynamics is leveraged on the order structure of the interval: each point separates the interval in two components. This feature does not exist for the plane and has to be replaced by a different separation property.
At any point $x$, one can consider its \emph{stable set}, i.e. the set of points whose iterates get closer to the forward orbit of $x$:
$$W^s(x)=\{y: {\rm dist}(f^n(x), f^n(y))\to 0\}.$$
For instance $x$ can belong to a periodic orbit which attracts all the points in a neighborhood: in this case $W^s(x)$ contains
a neighborhood of the orbit and one says that $x$ is a sink.

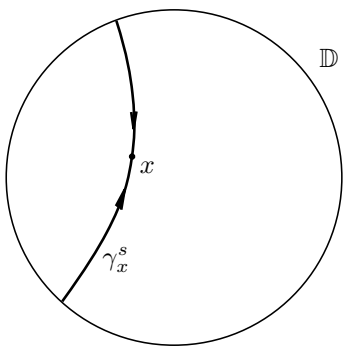
\begin{figure}[ht]
\begin{center}
\includegraphics[width=4.5cm]{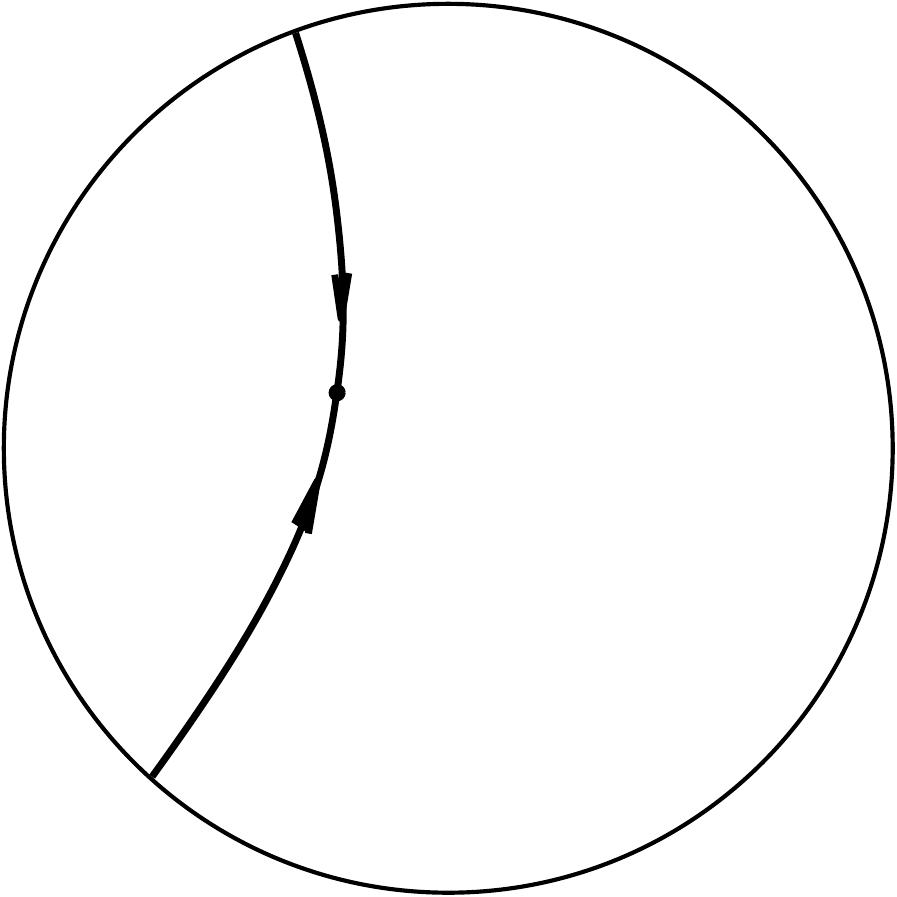}
\end{center}
\begin{picture}(0,0)
\put(-13,80){$x$}
\put(-27,45){$\gamma^s_x$}
\put(55,120){$\DD$}
\end{picture}
\caption{A stable manifold which separates the disc.
\label{f.strong-dissipation}}
\end{figure}

Since the dynamics is dissipative, one expects that for ``most points" $x$ the set $W^s(x)$ is non-empty, and indeed
using results from ergodic theory, one can show that unless $x$ is a sink, the stable set is an embedded $1$-dimensional submanifold, called the \emph{stable manifold} of $x$. Since the dynamics acts on the disc $\mathbb{D}$, we say that $W^s(x)$ separates, when it contains a curve $\gamma^s_x$ through $x$ whose endpoints belong to the boundary of the disc, so that $\DD\setminus \gamma^s_x$ has two connected components.
In this way, we can introduce the  following definition (see figure~\ref{f.strong-dissipation}):
\medskip

\noindent
\textbf{Definition.}
\emph{A diffeomorphism which sends the disc into its interior and contracts the volume is \emph{mildly dissipative} if for any invariant\footnote{The invariance of the measure by $f$ means that $f^*\mu=\mu$.} probability measure $\mu$
and for $\mu$-almost every point $x$,
\begin{itemize}
\item either $x$ is an attracting periodic point (a sink),
\item or through $x$ there exists a curve that is contained in the stable set of $x$ and that separates the disc.  
\end{itemize}}

It turns out that this class contains open sets of maps that are sufficiently close to maps on the interval
(hence contains strongly dissipative systems), but is wider: using tools from complex analysis, one can show
that it also contains all H\'enon maps whose Jacobian have modulus less than $1/4$.
For this reason these systems are called mildly dissipative.
{\color{black} Note that one can build examples of dissipative diffeomorphisms on the disc that are not mildly dissipative,
but these systems are maybe exceptional: we do not know if mild dissipation holds generically.}

\subsection*{One-dimensional reduction}
{\color{black} Although the existence of stable curve only occurs on a measurable subset, it allows to induce dynamical partitions of the system.}
Assuming  the mild dissipation property (stable manifolds separate the disc) one can prove that the dynamics of a mildly dissipative diffeomorphisms of the disc can be reduced to a continuous non-invertible map acting on a real tree (a simply connected and path connected metric space):
\medskip

\noindent
\textbf{Property.}
\emph{Given a smooth mildly dissipative diffeomorphisms f of the disc $\mathbb{D}$, there exist
a continuous map $h$ on a real tree $X$ and a projection $\Pi\colon \DD\to X$ such that:
\begin{itemize}
\item $f$ and $h$ are semi-conjugated: $\Pi\circ f=h\circ \Pi$,
\item any two $f$-invariant probability measures $\mu,\nu$ with no atoms and mutually singular
project on different measures $\Pi_*(\mu),\Pi_*(\nu)$.
\end{itemize}}
\smallskip

The second item says that the projection does not collapse too much the dynamics.
\medskip

\begin{figure}[ht]
\begin{center}
\includegraphics[width=6cm]{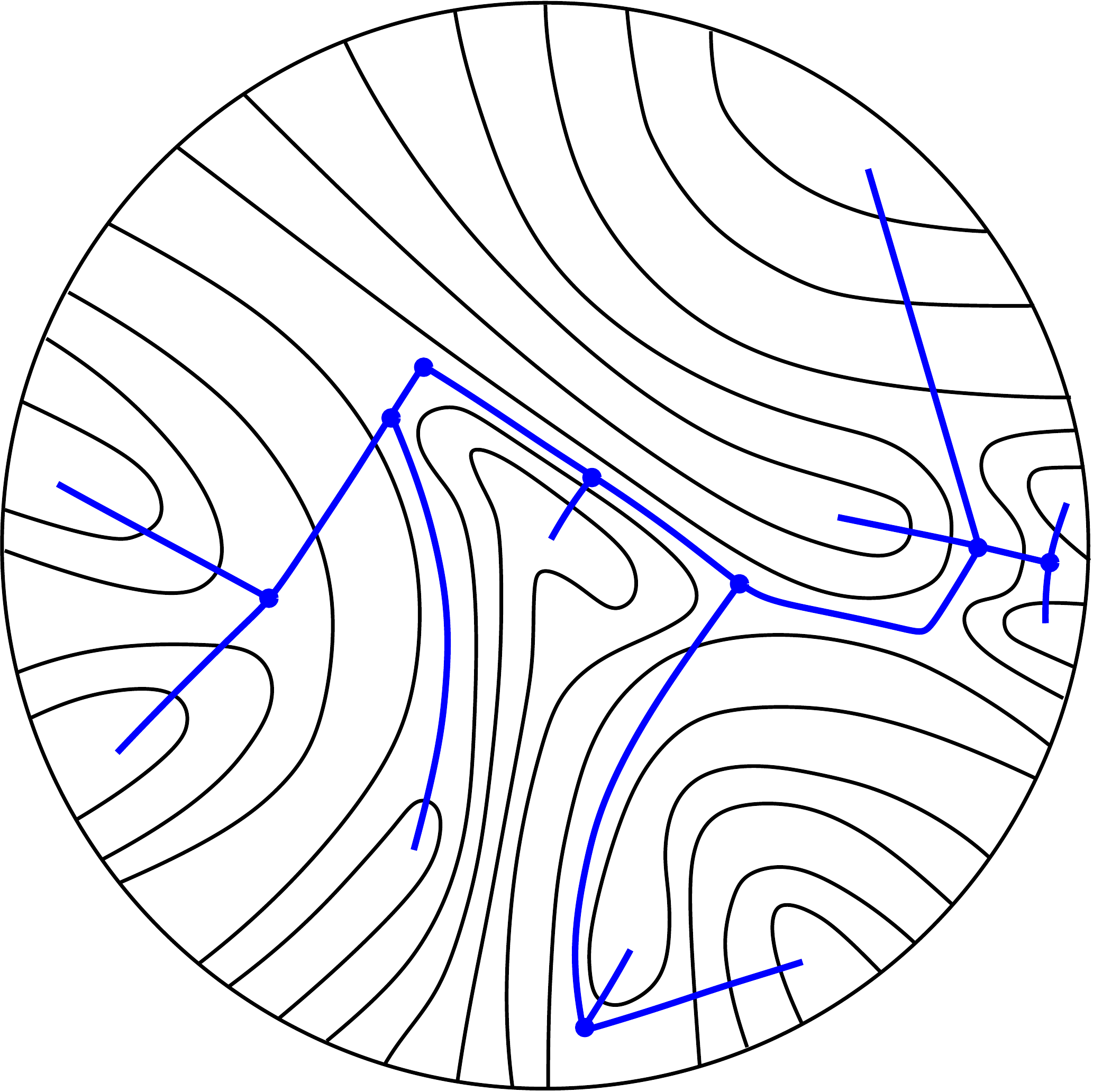}
\end{center}
\caption{ The one-dimensional structure associated to the family of stable manifolds. \label{f.lamination}}
\end{figure}

Reducing a system to a lower-dimensional one is a frequent strategy in dynamics.
In our setting the key idea behind is that the space of leaves of foliations in the plane generates a one-dimensional structure.
The mild dissipation provides a large collection of stable manifolds that are disjoint separating curves. It is well-known that the dual object to a planar lamination is a tree: the idea behind the proof of the previous property is to quotient the disc along these stable curves, see figure~\ref{f.lamination}.
That property suggests that the one-dimensional order structure re-emerges from the mild dissipation and makes possible to envision results with a ``one-dimensional flavor''.

\subsection*{Periodic approximation in the disc}
Another concise result that highlights the richness of the mildly dissipative class is the following:
\medskip

\noindent
\textbf{Property.}
\emph{For mildly dissipative diffeomorphisms of the disc, the closure of the set of periodic points contains the support of any invariant probability measure.}
\medskip

The proof uses an essential idea of one-dimensional dynamic that can be transposed to mildly dissipative diffeomorphisms of the disc: recurrence of non-periodic points forces to reverse the orientation on the projected tree $X$ and this implies the existence of a periodic point. 
It goes along the following lines.

\begin{figure}[ht]
\label{closing}
\begin{center}
\includegraphics[width=5.5cm]{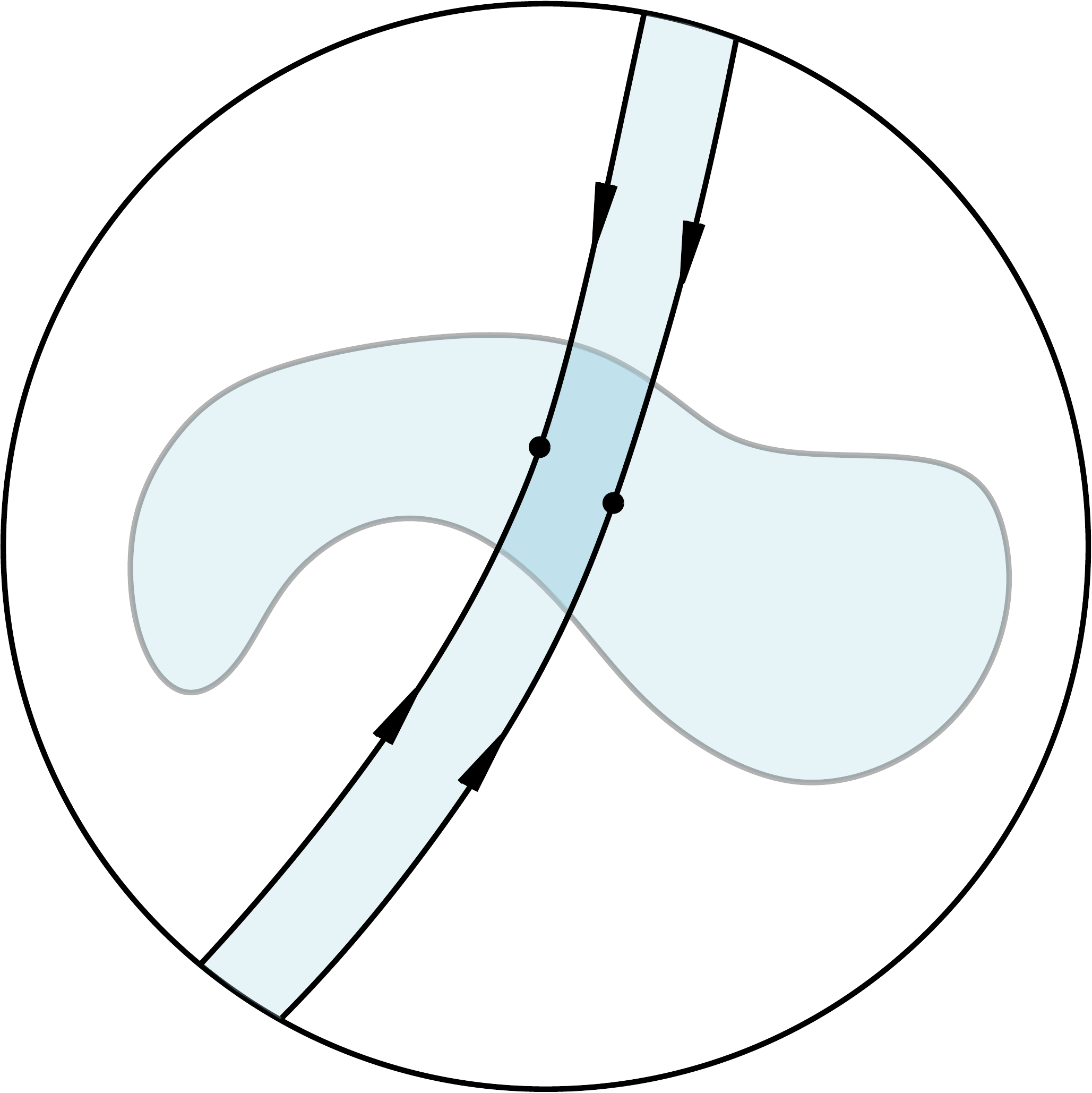}
\begin{picture}(0,0)
\put(-120,85){$D^-$}
\put(-45,67){$D^+$}
\put(-69.5,85){\small $x_1$}
\put(-93.5,93){\small $x_0$}
\put(-77.5,96){\small $R$}
\put(-110,30){\small $S$}
\put(-145,135){$\DD$}
\end{picture}
\caption{Why periodic points are dense. \label{f.closing}}
\end{center}\end{figure}
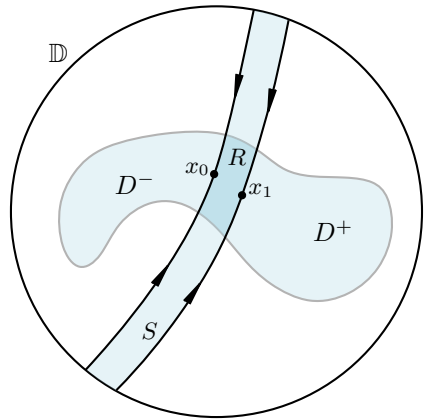

Let us consider an $f$-invariant probability measure $\mu$ with no atom.
Poincar\'e recurrence theorem asserts that, in a restriction to a full measure set, all points are recurrent.
Let us fix some point $x_0$ in that set. We have to prove that any neighborhood of $x_0$ contains a periodic point.
By recurrence, there exists $x_1=f^n(x_0)$ close to $x_0$ such that the stable curves $\gamma^s_{x_0}$, $\gamma^s_{x_1}$ are close and bound a thin strip $S$;
the intersection of the strip with a large iterate of the disc, $D:=f^m(\DD)$
defines a box $R$ with small diameter as in figure~\ref{f.closing}, such that $D\setminus R$ has two connected components $D^-$ and $D^+$.

Let $g:=f^n$ and let $k$ be the first positive integer such that $g^k(x_0)\in D^+$ and $g^{k+1}(x_0)\notin D^+$
(this exists since $x_1=g(x_0)$ does not belong to $D^+$ and $x_0$ is recurrent).
Similarly as in the one-dimensional case, we consider a continuous map $h\colon D\to R$
such that
$$h|_R=\operatorname{Id},\quad
h(D^-)= \gamma^s_{x_0}\cap R, \text{ and }h(D^+)= \gamma^s_{x_1}\cap R.$$
In particular, $h\circ g^k$ sends $R$ into itself and therefore has a fixed point $p$ in $R$.

Since $g^k(x_0)\in D^+$ and since stable curves do not cross, $g^k(\gamma^s_{x_0})\subset D^+$;
similarly $g^k(\gamma^s_{x_1})\subset D^-$. Consequently,
$h\circ g^k$ has no fixed point in $R\cap (\gamma^s_{x_0}\cup \gamma^s_{x_1})$ and by definition of $h$,
one deduces $h\circ g^k(p)=g^k(p)=p$. Hence $f$ has a periodic point in $R$ which is arbitrarily close to $x_0$, as required.

\section*{Mildly dissipative dynamics with zero entropy:\\ I-- prototype models}
We have described two classes of examples of mildly dissipative diffeomorphisms:
Morse-Smale systems, which belong to the interior of the set of dynamics with zero entropy,
and the examples pictured on the figure~\ref{f.odometer} which belong to its boundary.
We now present topological models with zero entropy and unbounded periods,
that can be built through a sequences of ``surgeries and pasting'' of two elementary Morse-Smale systems.

\subsection*{ Prototype models.}
\begin{figure}
\begin{center}
\includegraphics[width=7.8cm,height=3.3cm, angle=0]{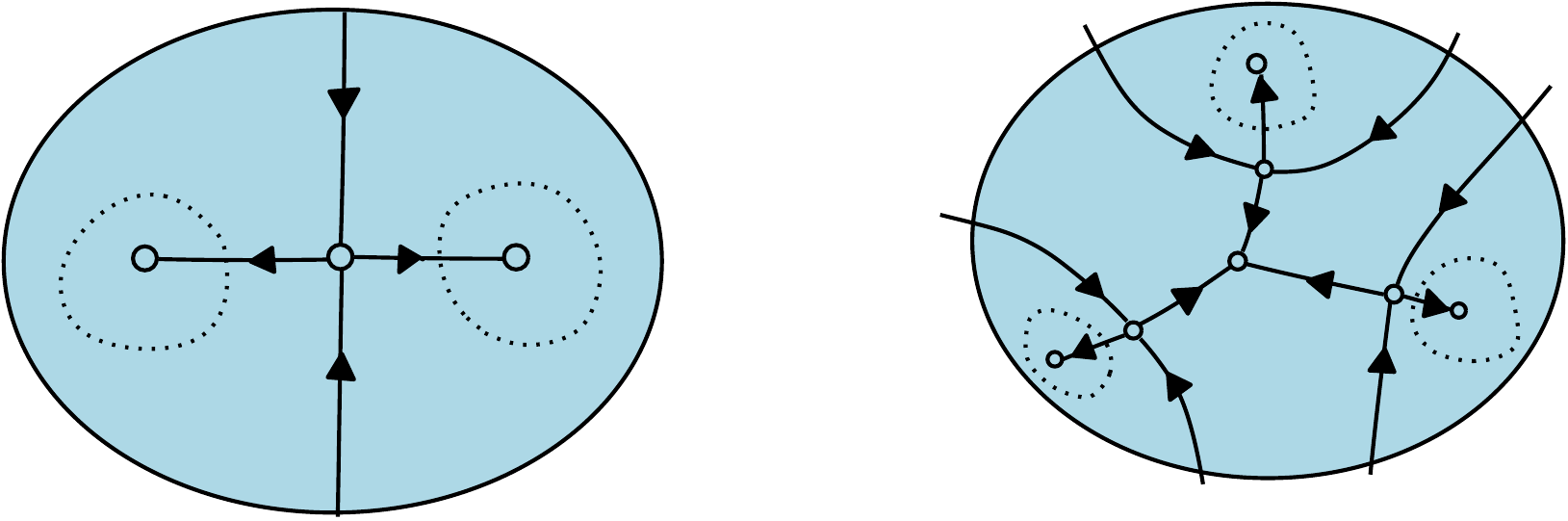}
\end{center}
\caption{The diffeomorphisms $f_0$ (left) and $f_1$ (right). The attracting domains are depicted by a dash boundary.\label{MSpdf}}
\end{figure}

Let us first introduce two Mor\-se-Smale dissipative diffeomorphisms of the disc $f_0, f_1$
that we describe below and depicted in figure~\ref{MSpdf}.
The limit set of $f_0$ is the union of a fixed saddle whose unstable branches are interchanged and of an attracting orbit of period two that revolves around the fixed point.
The limit set of $f_1$  is the union of a fixed attracting periodic point, a saddle of period three revolving around the fixed point
and an attracting periodic orbit (also of period three); each saddle has an unstable branch anchored at the fixed point and an unstable branch contained in the attracting domain of the $3$-periodic sink. Both diffeomorphisms are depicted in figure~\ref{MSpdf}.
In both situations, one says that the saddle periodic orbit is stabilized: either it is a fixed point, or its unstable manifold intersects
the basin of a fixed sink.

\subsubsection*{An inductive construction}
Given any sequence $(k_i)$ in $\{0,1\}^\NN$, one can build a sequence of dissipative diffeomorphisms
$(g_i)$, with exactly one sink of period $\tau_i:=\Pi_{j=1}^i(2+k_j)$ whose basin is a disc $D_i$.
It is obtained inductively from the diffeomorphism $g_{i-1}$ by ``pasting" the diffeomorphism $f_{k_i}$
in the basin of the sink of $g_{i-1}$, so that the return map $g_i^{\tau_{i-1}}|_{D_i}$ is conjugated to $f_{k_i}$.
In that way, $g_i$ has a nested sequence of discs $D_0\supset D_1\supset\dots\supset D_i$
that are renormalization domains of periods $\tau_0,\dots,\tau_i$.
Each diffeomorphism $g_i$ is Morse-Smale; moreover the construction can be done in a such way that the sequence $(g_i)$ converges to a homeomorphism.

\subsubsection*{Properties of the  limit system} The homeomorphism that is obtained as limit of the sequences $(g_i)$ verifies that:
\begin{itemize}
\item the dynamics is ``infinite  renormalizable'' in the sense that there is a nested sequences of renormalization domains with increasing periods; 
\item  the limit set (i.e. the set of points that belongs to renormalization domains with arbitrarily large period) is a Cantor set whose dynamics is an odometer (as introduced at the beginning of this text,
but its sequence of periods $(\tau_i)$ may not be equal to the sequence $(2^i)$).
\end{itemize}

We want to make some remarks: (i) The construction shows that there exist homeomorphisms with vanishing entropy and with periodic points whose period is not $2^n$.
(ii) The sequence can converge to a smooth  mildly dissipative diffeomorphism if $k_i=0$ for $i$ large.
(iii) The previous construction can be performed by gluing together more elementary diffeomorphisms $f_k$:
the period of their saddles and of their non-fixed sink may be larger; one can also consider more complicate Morse-Smale systems $f_0,f_1$.

\subsection*{Are the prototype models the typical ones?}
One can ask if the properties displayed by the prototype models are also satisfied by mildly dissipative diffeomorphisms $f$ of the disc with zero entropy.  More precisely:

\begin{itemize}
 \item What can be the periods of a nested sequence of attracting domains?
 \item When $f$ belongs to the interior of the set of systems with entropy zero,
do its periodic points have bounded periods?
 \item When $f$ belongs to the boundary of the set of systems with entropy zero,
is it infinitely renormalizable? is any limit set either an odometer Cantor set or a periodic orbit?
\end{itemize}

\subsection*{The strongly dissipative case}
These questions can be tested on strongly dissipative diffeomorphisms.
In fact, they have been answered by de Carvalho, Lyubich and Martens~\cite{dCLM} for H\'enon-like mappings of the form $F(x, y) = (f(x) + \eps(x, y), x)$ where $f$ is a unimodal map of the interval with a quadratic turning point and $\eps$ is  a real-valued map from the square to $\RR$ with a small size. In this work, they construct a period-doubling renormalization operator which extends the renormalization operator
introduced for unimodal maps (figure~\ref{f.renormalization-global}) and they show that (for sufficiently small $\epsilon$)
the properties carry over to this case. Namely, the renormalization operator admits a unique fixed point (which actually coincides with the fixed point of the renormalization on the interval): it is hyperbolic (with a one-codimensional stable manifold)
and the periods of its renormalization domains are $2^n$ for all $n\geq 0$. 

\section*{Mildly dissipative dynamics with zero entropy:\\ II-- the general case}
In the general case, the notion of turning point does not exist anymore and the map may be far from one-dimensional endomorphisms.  In fact, it is not difficult to construct mildly dissipative diffeomorphisms with zero entropy which are not close to an interval map and even have periodic points with periods that are not a power of two (see the prototype construction) and so the renormalization scheme developed for H\'enon-like maps with very small Jacobian cannot be applied directly.

In what follows, we are going to state the results that we have obtained with Charles Tresser~\cite{CPT} and at the end, we explain some of the main ideas of the proof.

\subsection*{Renormalizable dynamics.}
As in dimension $1$, the re\-nor\-malization is an essential tool for secribing the transition to chaos. Let us define that notion for surface diffeomorphisms.
\medskip

\noindent
{\bf Definition.}
 A diffeomorphism $f$ of the disc is \emph{renormalizable} if
there exist a compact set $D\subset \mathbb{D}$ homeomorphic to the unit disc
and an integer $\tau>1$ such that
$f^i(D)\cap D=\emptyset$ for each $1\leq i<\tau$ and $f^\tau(D)\subset D$.
One says that $D$ is a \emph{renormalization domain of period $\tau$}.
\medskip

Based on that definition, one gets a dichotomy:

\begin{theorem}\label{t.theoremA}
For any mildly dissipative diffeomorphism $f$ of the disc whose entropy vanishes,
\begin{itemize}
\item either $f$ is renormalizable,
\item or any forward orbit converges to a fixed point.
\end{itemize}
\end{theorem}

Morse-Smale diffeomorphisms (whose non-wandering dynamics is a finite set of hyperbolic periodic points)
are certainly not infinitely renormalizable.
It is natural to generalize this class of diffeomorphisms
in order to allow bifurcations of periodic orbits.
\medskip

\noindent
{\bf Definition.}
A system is \emph{generalized Morse-Smale} if:
\begin{itemize}
\item the limit set of any forward orbit is a periodic orbit,
\item the limit set of any backward orbit that is contained in $\mathbb{D}$ is a periodic orbit,
\item the set of periods over all periodic orbits is finite.
\end{itemize}
(Contrary to classical Morse-Smale systems, there may exist infinitely many periodic points with a same period.)
\medskip

Clearly these diffeomorphisms have zero entropy. Moreover, the set of mildly dissipative generalized Morse-Smale diffeomorphisms of the disc
is $C^1$ open. A stronger version of theorem~\ref{t.theoremA} states that in the renormalizable case there exist finitely many disjoint renormalization domains  such that the limit set of any forward orbit contained in their complement is a fixed point. That version  implies:
\medskip

\noindent
{\bf Corollary.}
\emph{A mildly dissipative diffeomorphism of the disc with zero entropy is
\begin{itemize}
\item either infinitely renormalizable,
\item or generalized Morse-Smale.
\end{itemize}}

\subsection*{ Boundary of zero entropy.}
From the previous theorem and the fact that generalized Morse-Smale diffeomorphisms are in the interior of the set of systems with zero entropy, one can characterize the dynamics in the boundary of zero entropy:
\medskip

\noindent
{\bf Corollary.}
\emph{A mildly dissipative diffeomorphism of the disc in the boundary of zero entropy
is infinitely renormalizable.}
\medskip

One can now wonder, after these results, if one can get a complete characterization of the limit sets of these systems.
The following result extends the property of interval maps.
\medskip

\noindent
{\bf Corollary.}
\emph{Let $f$ be a mildly dissipative diffeomorphism of the disc with zero entropy.
Then the limit set of any orbit  is:
\begin{itemize}
\item either a periodic orbit,
\item or a \emph{generalized odometer}:
there exists an odometer $h$ on the Cantor set $\mathcal{K}$ and a
continuous subjective map $\pi\colon \mathcal{C}\to \mathcal{K}$ such that
$\pi\circ f|_\mathcal{C}=h\circ \pi$.
Moreover $\pi$ is ``essentially" one-to-one.
\end{itemize}}
\medskip

Figure~\ref{f.odometer} represents the second case.
\smallskip

\subsection*{Set of periods.}
One cannot expect that Sarkovskii's property stated above for interval maps extends identically in the disc.
Indeed the prototype models show that any finite set of integers can appear inside the set of periods of a
mildly dissipative diffeomorphism having zero entropy.
But a constraint appears, when one considers periodic orbits with sufficiently large period:

\begin{theorem}\label{t.period}
If $f$ is an infinitely renormalizable mildly dissipative diffeomorphism of the disc with zero entropy,
then there exist an open set $W$ and $m\geq 1$ such that:
\begin{itemize}
\item $W$ is a finite disjoint union of renormalization domains whose period divides $m$
(possibly several orbits of domains),
\item the periods of points in $\DD\setminus W$ are bounded by $m$,
\item any renormalization domain $D\subset W$ of $f^m$ has period
of the form $2^k$:
it is associated to a nested sequence of renormalization domains
$D=D_k\subset\dots \subset D_1\subset  W$ of $f^m$ with period $2^{k},\dots, 2$.
\end{itemize}
\end{theorem}
In other words, the period of a renormalization domain is eventually a power of $2$: after replacing $f$ by an iterate, the period of all the renormalization domains are powers of $2$. This implies the announced property on periods:
\medskip

\noindent
{\bf Corollary.}
\emph{For any mildly dissipative diffeomorphism $f$ of the disc with zero entropy,
there exist two finite families of integers $\{n_1,\dots,n_k\}$
and $\{m_1,\dots,m_\ell\}$
such that the set of periods of the periodic orbits of $f$ coincides with}
\begin{equation*}\label{e.period}
\operatorname{Per}(f)=\{n_1,\dots,n_k\}\cup \left\{m_i.2^k, \; 1\leq i\leq \ell \text{ and } k\in \mathbb{N}\right\}.
\end{equation*}

In particular, this proves Tresser's conjecture in the case of mildly dissipative dynamics of the disc. 
\smallskip

\subsection*{H\'enon maps}
The previous results can be applied to the H\'enon family for all parameters provided that the Jacobian is smaller than $1/4$ (this requires some adaptation in order to reduce it to a map sending the disc into its interior).
More precisely, when the entropy vanishes, any forward (resp. backward) orbit in $\RR^2$
has exactly one of the following behavior:
\begin{itemize}
\item it escapes to infinity, i.e. leaves compact sets;
\item it converges to a periodic orbit;
\item it accumulates to a generalized odometer.
\end{itemize}

\section*{Mildly dissipative dynamics with zero entropy:\\ III-- sketch of the proofs}
The approach for the general case cannot use the interval ordering and is based instead on the structure of the set of periodic points: the unstable branches of the saddle periodic points serve as a skeleton of the dynamics that allows to construct the renormalization domains. 
We fist explain this strategy on the prototype examples introduced before.

\subsection*{Dynamical features of the prototype examples.} 
\begin{figure}[h]
\begin{center}
\includegraphics[width=7cm, height=4cm, angle=0]{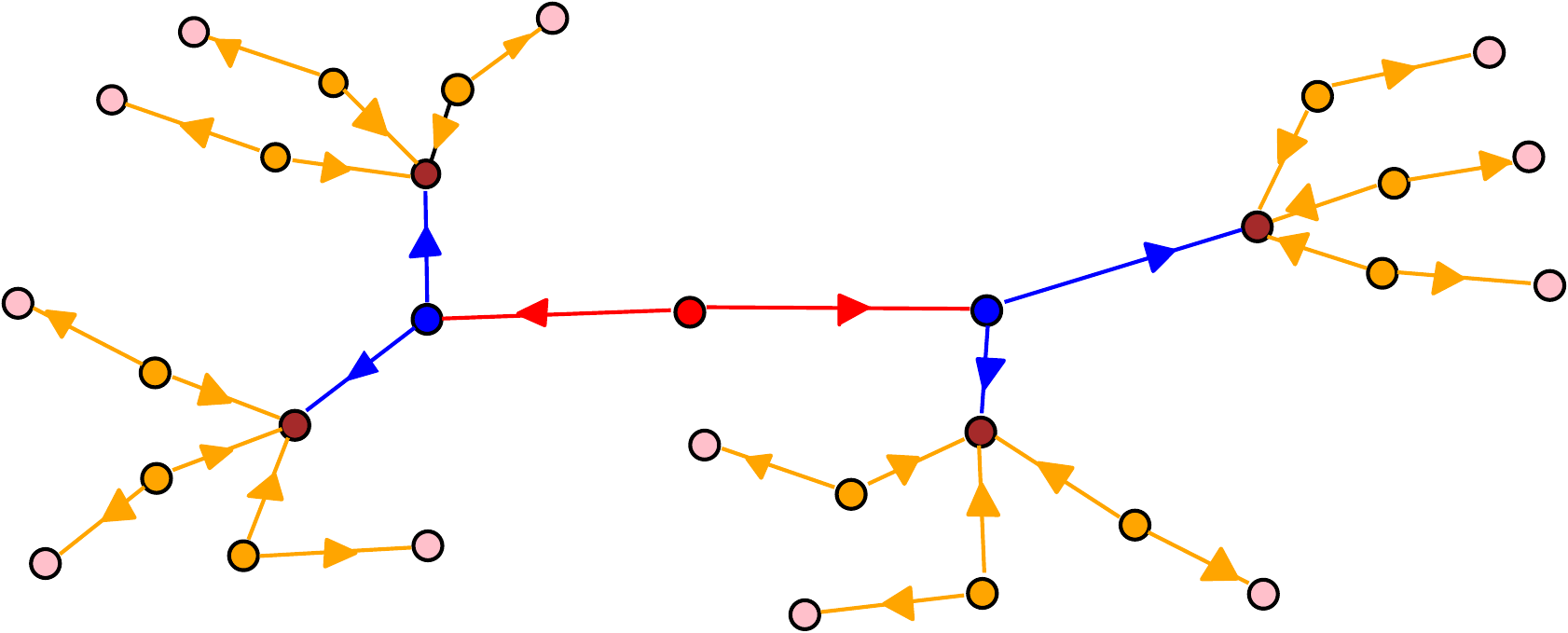}
\end{center}
\caption{Chain of periodic points associated to the diffeomorphism $g$ obtained by pasting successively $f_0,f_0,f_1$.
It contains:
one saddle fixed point (red), a two-periodic saddle orbit (blue) whose unstable branches are exchanged by $g^2$, a four-periodic attracting orbit (brown),
and twelve-periodic saddle (orange) and attracting (pink) orbits.
The arrows indicate if periodic points are saddles or sinks (for sinks, all arrows are pointing in).
\label{treepdf}}
\end{figure}

Let us consider a prototype diffeomorphism $g_i$ obtained after pasting a finite number Morse-Smale diffeomorphisms  $f_{k_0}, f_{k_{1}}, f_{k_i}$. The unstable branches of the saddles connect the periodic points and define a tree structure that we call {\it chain}, see figure~\ref{treepdf}.

\begin{figure}[h]
\begin{center}
\includegraphics[width=8cm,angle=0]{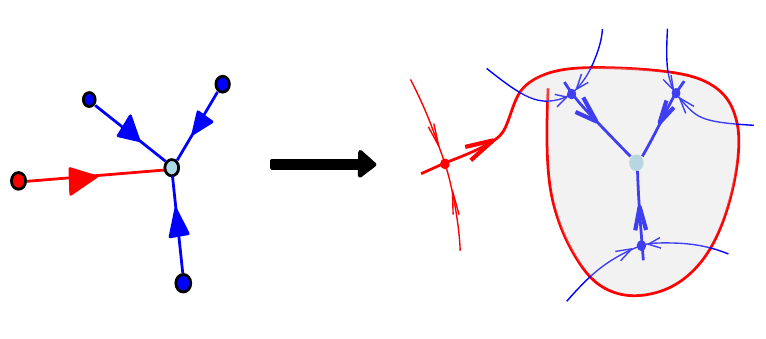}
\put(-230,53){$p$}
\put(-172,48){$q$}
\put(-172,8){$w$}
\put(-105,53){$p$}
\put(-35,48){$q$}
\put(-35,20){$w$}
\bigskip
\bigskip

\includegraphics[width=8cm,angle=0]{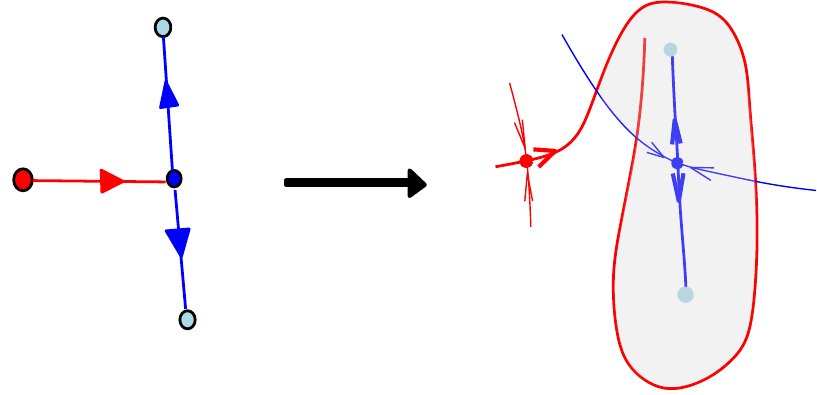}
\put(-228,65){$p$}
\put(-174,58){$q$}
\put(-90,70){$p$}
\put(-35,55){$q$}
\end{center}
\caption{Examples of Pixton discs.  \label{Pixtpdf}}
\end{figure}
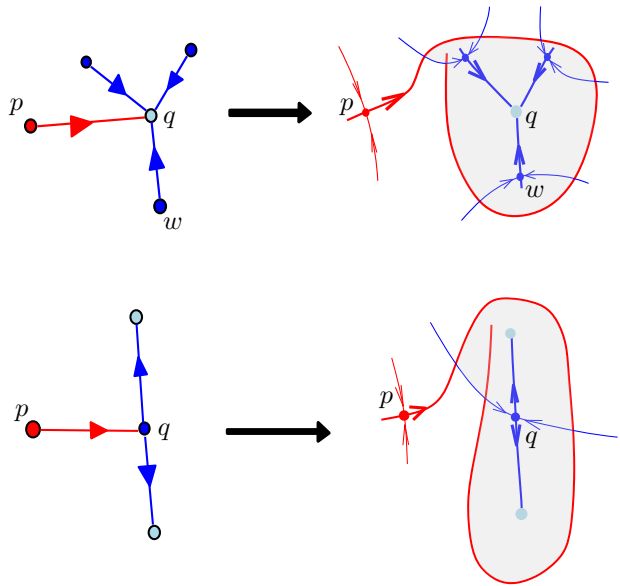

From each saddle $p$ points out at least one arrow, which lands at a point $q$ with the same or double period.
Two cases may occur (see figure~\ref{treepdf}):
\begin{itemize}
\item either $q$ is a an attracting periodic point,
\item or $q$ is a saddle whose unstable branches are exchanged by some iterate of $f$.
\end{itemize}
That observation allows to reconstruct the renormalization domains of the prototype example $f$, see figure~\ref{Pixtpdf}.

In the first case (above on the figure~\ref{Pixtpdf}), the unstable manifold of $p$ accumulates on the sink $q$ which anchors a revolving saddle $w$ with larger period (period three in the figure); this implies that  the unstable branch of $p$ has to cross the stable manifolds of the iterates of $w$.
One then defines a disc which contains $q,w$, is bounded by a piece of the unstable branch of $p$ and a piece of the stable manifolds of the saddle $w$, and which is mapped into itself by some iterate of $f$.

In the second case (below on the figure~\ref{Pixtpdf}), the unstable manifold of $p$ accumulates on the saddle $q$ (with the same period) whose unstable branches are exchanged by the dynamics and accumulate on a sink of double period. This implies that the unstable branch of $p$ has to cross both stable branches of $q$. Again, a piece of the unstable branch of $p$ and of the stable manifolds of the saddle $q$, defines a disc
which is mapped into itself by some iterate, contains $q$, but not $p$.

This construction leads us to introduce the following:
\smallskip

\noindent{\bf Definition.}
A Jordan domain $D$ is a {\em Pixton disc}\footnote{Pixton introduced a similar notion in order to study planar homoclinic orbits.}
of period $\tau$
if its boundary decomposes in two parts: one subset of $W^u(p)$ and a closed set whose forward iterates by $f^\tau$ are all contained in the interior of $D$.
\medskip

An \emph{trapped disc} for $f^\tau$ (i.e. a disc mapped by $f^\tau$ into its interior) is a particular example of a Pixton disc.
\medskip

\subsection*{How to work out the general case.}
The strategy in the general case goes along the next steps which will be detailed in the following paragraphs. 
\begin{itemize}
\item[i)] {\it Chains.}
As for the prototype models built previously, the set of fixed points and their unstable branches forms a connected set which has a tree structure. Considering also iterates $f^m$, one gets chains between periodic points whose period divides $m$:
the periodic points of larger period are connected to the ones of lower period and ``revolve'' around them. 
 
\item[ii)]{\it The case where all the periodic points are  fixed.} We then prove that any limit set is a fixed point.
 
\item[iii)]{\it Construction of Pixton discs.}
When there are periodic points that are not fixed, one builds Pixton discs
which contains all periodic points of higher period and are good candidates to be renormalization domains.
 
\item[iv)]{\it Renormalization domains.} Once the Pixton discs are constructed, we prove that the ``maximal ones" are renormalization domains.
 
%
\item[v)]{\it Eventual period two.}
At last one concludes that after several renormalizations, the new renormalization periods are all equal to $2$.
%
%
\end{itemize}

\subsection*{Chains of periodic points}
The key ingredient to obtain the tree structure is to check that there is no cycle between fixed (or periodic) points:
\medskip

\noindent
{\bf Property.}
\emph{There is no sequence of saddle fixed points $\{p_1,\dots, p_n\}$ such that the unstable manifold $W^u(p_i)$ accumulates on $p_{i+1}$ and
$W^u(p_n)$ accumulates on $p_1$.}
\medskip

This property generalizes Smale's theorem mentioned in the first section:
a cycle would force a situation close to what is depicted in figure~\ref{f.homoclinic}, which would give positive entropy.
\smallskip

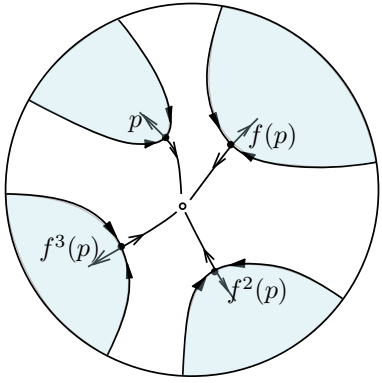
\begin{figure}
\begin{center}
\includegraphics[width=5cm,angle=0]{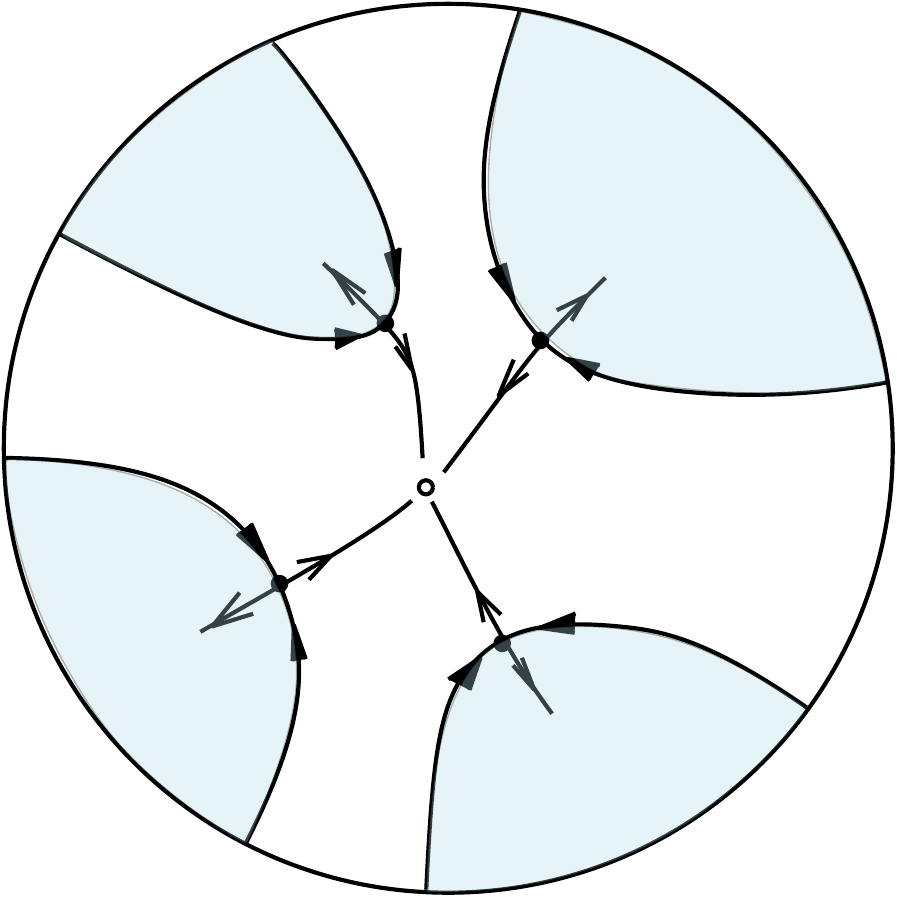}
\put(-58,30){$f^2(p)$}
\put(-50,89){$f(p)$}
\put(-95,95){$p$}
\put(-129,46){$f^3(p)$}
\end{center}
\caption{A stabilized periodic orbit of period $4$ and its $4$ decorated regions bounded by the stable manifolds.
The periodic orbit is stabilized by the fixed point in the middle. \label{f.decoration}}
\end{figure}

In chains, a  special role is played by \emph{stabilized points}: these are saddles that either are fixed and whose unstable branches are exchanged by $f$,
or are not fixed but whose unstable manifold is anchored by a fixed point.
The stable manifolds of the stabilized points bound domains called {\it decorated regions} (see figure~\ref{f.decoration}). These regions are pairwise disjoints: otherwise, using that each iterate of the decorated orbit has an unstable branch which accumulates on a stabilizing fixed point, it would imply that an unstable manifold crosses the stable manifold of another iterate, creating a homoclinic intersection and therefore contradicting the fact that the entropy vanishes.
 
Moreover, the decorating regions contain all the periodic points of larger periods: otherwise, it would again force a homoclinic intersection.
One can thus decompose the set of periodic points as:
\begin{itemize}
\item stabilizing fixed points,
\item stabilized periodic orbits,
\item periodic orbits contained in decorated regions.
\end{itemize}

\subsection*{ The case where all the periodic points are fixed.}
In this setting, the property of periodic approximation in the disc implies that any probability invariant measure is supported on the set of fixed points.
Hence, the limit set of any forward orbit contains a fixed point. If it is not a singleton, the forward orbit also accumulates on
points in unstable branches of fixed points, so that its limit set contains contains a cycle of fixed points. This would contradicts the no cycle property
stated before.

\subsection*{Construction of Pixton discs.}
\begin{figure}
\begin{center}
\includegraphics[width=8cm,angle=0]{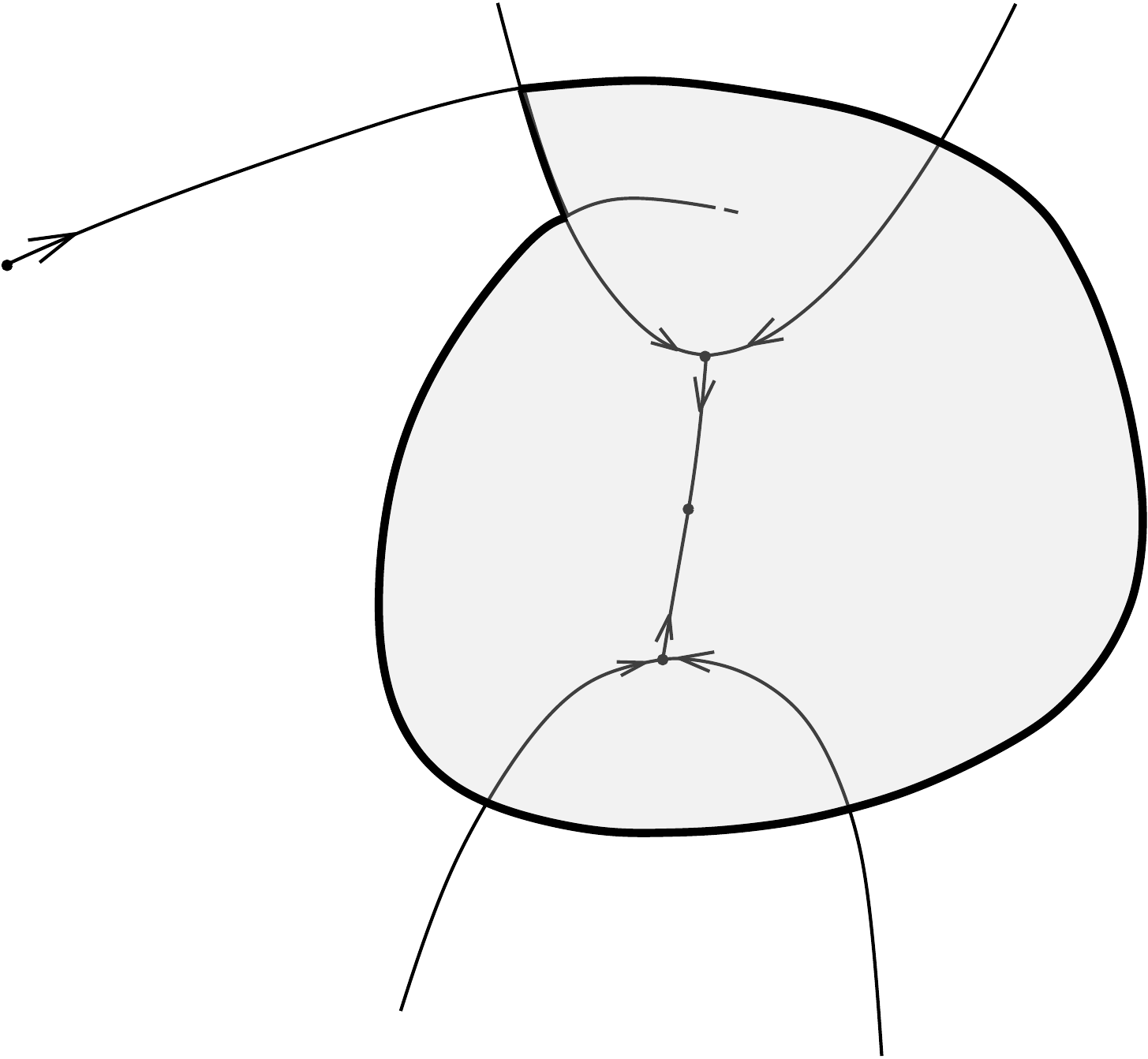}
\put(-178,168){\small $\Gamma$}
\put(-88,110){\small $q$}
\put(-220,150){\small $p$}
\put(-83,133){\small $w$}
\put(-92,83){\small $f(w)$}
\put(-35,80){\small $D$}
\end{center}
\caption{Construction of a Pixton disc: $\Gamma$ is fixed and $w$ has period $2$.\label{f.pixton}}
\end{figure}
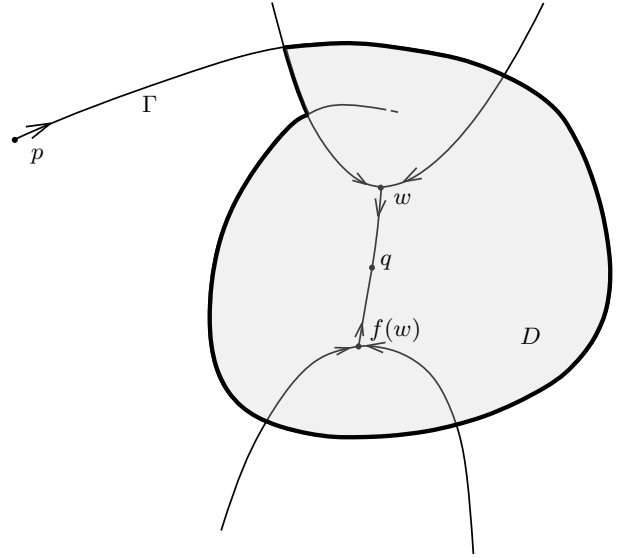

To each unstable branch $\Gamma\subset W^u(p)$, fixed by an iterate $f^\tau$, we build a Pixton disc $D_\Gamma$ for $f^\tau$ that contains the accumulation set of $\Gamma$, in a similar way as we did for the prototype examples:
if $w$ is a saddle point accumulated by $\Gamma$, one considers a disc $D$ bounded by an arc in $\Gamma$ and an arc in the stable manifolds of $w$; this disc contains all the periodic points of deeper level and connected to $w$ in the chain structure, see figure~\ref{f.pixton}.
The Pixton disc $D_\Gamma$ is obtained as the union of such discs for different choices of $w$.

\subsection*{Renormalization  domains.}  To prove that the Pixton disc $D_\Gamma$ is actually a renormalization domain,
one has to prove that the iterates of $\Gamma\cap D_\Gamma$ (in figure~\ref{f.pixton}) remain contained in the disc; if a piece of $\Gamma$ escapes from $D_\Gamma$ under forward iterations, a strong version of the property of periodic approximation
implies that there are periodic points outside $D_\Gamma$
which are accumulated by $\Gamma$, a contradiction since $D_\Gamma$
contains all the periodic points that belong to the accumulation set of $\Gamma$.


\subsection*{Eventual period two.}
The previous steps build the renormalization (theorem A).
A large number of renormalizations reduces the study to a small neighborhood $W$ of the union of the generalized odometers.
We then have to show that the period of all the further renormalizations is equal two (theorem \ref{t.period}).

We first observe that for the saddle orbits contained in $W$, a large proportion of the iterates have stable manifolds which vary continuously
for the $C^1$-topology. In particular, for a large proportion of points, the stable manifolds are ``parallel". This is consistent with the example of figure~\ref{f.odometer}, where the renormalization periods are $2$ at each step. 
However a renormalization period larger than two would provide more than two stable curves, based at iterates close, and which have to bend
away from each other  (see for instance the figure~\ref{f.decoration} where the period is $4$). This contradicts the fact that these curves
are $C^1$-close.

\section*{Dynamics in higher dimensions}
There is no such detailed description of the dynamics of systems with zero entropy
for general surface diffeomorphisms and on higher-dimensional manifolds.
However perturbative methods have been developed which allow to describe a $C^1$-dense open set of systems.
In particular, they imply the following dichotomy:
\medskip

\noindent
{\bf Theorem} \cite{PuSa,C}.
\emph{The union of the set of Morse-Smale diffeomorphisms and
of the set of diffeomorphisms having a transverse homoclinic intersection is
a $C^1$-dense open set of the space of diffeomorphisms.}
\medskip

As the diffeomorphisms with a transverse homoclinic intersection have positive entropy,
this result characterizes - inside a dense open set - the systems with zero entropy.
However the dynamics on the boundary of the set of systems with zero entropy are not understood.
And in higher topologies, almost nothing is known.
\bigskip

{\color{black}
\noindent
\emph{Acknowledgments. We are grateful to the anonymous referees for their numerous comments on this text which helpt to improve its presentation.}
}

\bigskip
\bigskip

\begin{tabular}{l l}
\emph{Sylvain Crovisier}
&
\emph{Enrique Pujals}
\\

LMO - CNRS
& Graduate Center \\
Univ. Paris-Saclay
& CUNY\\
Orsay, France
&  New York, USA
\end{tabular}
 
 \end{document}